\documentclass[opre]{informs3}
\makeatletter
\newcommand{\reflabel}[2]{#2\def\@currentlabel{#2}\label{#1}}
\makeatother

\newcommand{\E}{\mathbb E}
\renewcommand{\P}{\mathbb P}
\newcommand{\Q}{\mathbb Q}
\newcommand{\R}{\mathbb R}

\newcommand{\hX}{\widehat X}

\newcommand{\hZ}{\widehat Z}

\newcommand{\sF}{\mathscr F}
\newcommand{\cG}{\mathcal G}
\newcommand{\cK}{\mathcal K _c}
\newcommand{\cP}{\mathcal P}

\newcommand{\cL}{\mathcal L}
\newcommand{\cW}{\mathcal W}
\newcommand{\cX}{\mathcal X}
\newcommand{\cZ}{\mathcal Z}
\newcommand{\Fcomp}{\sF _{\hP}}

\newcommand{\RP}{\cL}
\newcommand{\Id}{\operatorname{Id}}
\newcommand{\Proj}{\operatorname{Proj}}
\newcommand{\Graph}{\operatorname{Graph}}

\newcommand{\cM}{\overline {\mathcal K} _c}
\newcommand{\tP}{\widetilde \P}
\newcommand{\tX}{\widetilde X}
\newcommand{\tx}{\tilde x}
\newcommand{\tf}{\tilde f}
\newcommand{\tc}{\tilde c}
\newcommand{\tcK}{\mathcal K _{\tilde c}}

\newcommand{\cPs}{\bar{\cP}}
\newcommand{\cSc}{\cS ^c}

\newcommand{\cLm}{\overline{\mathcal L}}
\newcommand{\cLJ}{\cL _{\mathrm J}}
\newcommand{\cLJm}{\cLm _{\mathrm J}}
\newcommand{\cLG}{\cL _{\mathrm G}}
\newcommand{\hQ}{\widehat {\mathbb Q}}
\newcommand{\halpha}{\widehat \alpha ^*}

\newcommand{\doublestar}{^{\substack{\vspace{-.2em}*\\{*}}}}
\newcommand{\MAD}{\operatorname{\mathbb MAD}}
\newcommand{\Ent}{\operatorname{\mathbb Ent}}

\newcommand{\Rninfty}{\R \cup \{-\infty\}}
\newcommand{\Rbar}{\bar \R}

\newcommand{\CVaRP}{\CVaR ^{\P}}
\newcommand{\CVaRhP}{\CVaR ^{\hP}}
\newcommand{\CVaRQ}{\CVaR ^{\Q}}
\newcommand{\CVaRhQ}{\CVaR ^{\hQ}}

\newcommand{\MADP}{\MAD ^{\P}}
\newcommand{\MADhP}{\MAD ^{\hP}}
\newcommand{\VarP}{\Var ^{\P}}
\newcommand{\VarhP}{\Var ^{\hP}}
\newcommand{\EntP}{\Ent ^{\P}}
\newcommand{\EnthP}{\Ent ^{\hP}}
\newcommand{\MADQ}{\MAD ^{\Q}}
\newcommand{\MADhQ}{\MAD ^{\hQ}}
\newcommand{\VarQ}{\Var ^{\Q}}
\newcommand{\VarhQ}{\Var ^{\hQ}}
\newcommand{\EntQ}{\Ent ^{\Q}}
\newcommand{\EnthQ}{\Ent ^{\hQ}}

\DeclareMathOperator*{\essssup}{ess\,sup}
\usepackage{multirow}
\SingleSpacedXI

\renewcommand{\theARTICLETOP}{}
\renewcommand{\theARTICLERULE}{\vspace{12pt}}
\renewcommand{\theRRHFirstLine}{}
\renewcommand{\theRRHSecondLine}{}
\renewcommand{\theLRHFirstLine}{}
\renewcommand{\theLRHSecondLine}{}
\renewcommand{\theECLRHFirstLine}{}%
\renewcommand{\theECLRHSecondLine}{}%
\renewcommand{\theECRRHFirstLine}{}%
\renewcommand{\theECRRHSecondLine}{}%
\renewcommand{\endproof}{\endTrivlist\addvspace{2ex}}
\makeatletter
\renewcommand\section{\@startsection {section}{1}{\z@}{-13pt plus -6pt minus -3pt}{4pt}%
  {\fs.13.15.\bfseries\RAGG}}%
\renewcommand\subsection{\@startsection{subsection}{2}{\z@}{-13pt plus -6pt minus -3pt}{4pt}%
  {\TEN\bfseries\RAGG}}%
\makeatother
\let\FigureCaptionFontStyle\undefined
\def\FigureCaptionFontStyle{\mdseries}

\usepackage[scaled=.98,sups,osf]{XCharter}%
\usepackage[utf8]{inputenc}
\usepackage[T1]{fontenc}
\usepackage[english]{babel}
\usepackage{microtype}
\usepackage[smallerops,vvarbb]{newtxmath}
\usepackage[varqu,varl]{inconsolata}%
\usepackage[scaled=0.95]{helvet}
\DeclareMathAlphabet{\mathsfit}{T1}{\sfdefault}{\mddefault}{\sldefault}
\SetMathAlphabet{\mathsfit}{bold}{T1}{\sfdefault}{\bfdefault}{\sldefault}
\DeclareMathAlphabet{\mathcal}{OMS}{cmsy}{m}{n}
\usepackage{mathtools,accents}
\usepackage{xifthen,enumitem,comment}
\setlist[enumerate,1]{label=\normalfont{(\Roman*)},leftmargin=2em}
\usepackage{xcolor}
\usepackage{etoolbox}
\makeatletter
\patchcmd{\env@cases}{1.2}{0.96}{}{}
\makeatother
\usepackage{algorithm,algpseudocode}
\usepackage[hypertexnames=false]{hyperref}
\hypersetup{%
  colorlinks=true,%
  linkcolor={blue!66!black},%
  linkbordercolor=red,%
  citecolor={blue!66!black},
}
\usepackage{booktabs,tikz,pgfplots,subfigure}
\usetikzlibrary{matrix,decorations.pathreplacing}

\ifx\argmax\undefined\DeclareMathOperator*{\argmax}{arg\,max}\fi
\ifx\argmin\undefined\DeclareMathOperator*{\argmin}{arg\,min}\fi

\newcommand*{\QED}{%
\leavevmode\unskip\penalty9999 \hbox{}\nobreak\hfill
    \quad\hbox{$\square$}%
}
\newcommand*{\QEG}{%
\leavevmode\unskip\penalty9999 \hbox{}\nobreak\hfill
    \quad\hbox{$\clubsuit$}%
}
\newcommand*{\QDEF}{%
\leavevmode\unskip\penalty9999 \hbox{}\nobreak\hfill
    \quad\hbox{$\diamondsuit$}%
}

\renewcommand{\Pr}{\mathbb{P}}
\renewcommand{\P}{\mathbb{P}}
\newcommand{\hP}{\widehat{\mathbb{P}}}
\newcommand{\hx}{\widehat{x}}

\providecommand{\E}{\mathbb{E}}
\providecommand{\R}{\mathbb{R}}
\providecommand{\CVaR}{\mathbb{C}\mathrm{V@R}}
\providecommand{\Var}{\mathbb{V}\mathrm{ar}}
\providecommand{\Q}{\mathbb{Q}}

\providecommand{\frakM}{\mathfrak{M}}
\providecommand{\cA}{\mathcal{A}}

\providecommand{\cG}{\mathcal{G}}

\providecommand{\cK}{\mathcal{K}}
\providecommand{\cL}{\mathcal{L}}
\providecommand{\cM}{\mathcal{M}}

\providecommand{\cP}{\mathcal{P}}

\providecommand{\cS}{\mathcal{S}}
\providecommand{\cT}{\mathcal{T}}
\providecommand{\cU}{\mathcal{U}}

\providecommand{\cW}{\mathcal{W}}
\providecommand{\cX}{\mathcal{X}}

\providecommand{\cZ}{\mathcal{Z}}

\providecommand{\Dirac}{\boldsymbol{\updelta}}

\providecommand{\supp}{\mathrm{supp}\,}

\usepackage[numbers]{natbib}

\theoremstyle{THkey}\newtheorem{IP}{IP}

\TheoremsNumberedThrough     %
\ECRepeatTheorems

\EquationsNumberedThrough    %

\pgfplotsset{compat=1.17} 
\begin{document}

\RUNAUTHOR{Zhang, Yang, and Gao}

\RUNTITLE{Duality for Wasserstein DRO}

\TITLE{A Short and General Duality Proof for\\ Wasserstein Distributionally Robust Optimization}

\ARTICLEAUTHORS{%
\AUTHOR{Luhao Zhang}
\AFF{Department of Industrial Engineering and Operations Research,
Columbia University\\ %
\EMAIL{lz2487@columbia.edu}} %

\AUTHOR{Jincheng Yang}
\AFF{Department of Mathematics,
University of Chicago\\ %
\EMAIL{jincheng@uchicago.edu}} %

\AUTHOR{Rui Gao}
\AFF{Department of Information, Risk and Operations Management, 
The Unversity of Texas at Austin\\ %
\EMAIL{rui.gao@mccombs.utexas.edu}} %
} %

\ABSTRACT{%
We present a general duality result for Wasserstein distributionally robust optimization that holds for any Kantorovich transport cost, measurable loss function, and nominal probability distribution. 
Assuming an interchangeability principle inherent in existing duality results, our proof only uses one-dimensional convex analysis. Furthermore, we demonstrate that the interchangeability principle holds if and only if certain measurable projection and weak measurable selection conditions are satisfied.
To illustrate the broader applicability of our approach, we provide a rigorous treatment of duality results in distributionally robust Markov decision processes and distributionally robust multistage stochastic programming.
Additionally, we extend our analysis to other problems such as infinity-Wasserstein distributionally robust optimization, risk-averse optimization, and globalized distributionally robust counterpart.
}%


\maketitle

\section{Introduction}

In this paper, we consider the following problem 
\begin{align}
    \label{problem}
    \RP (\rho) := \sup _{\P \in \cP (\cX)} \left\lbrace
        \E _{X \sim \P} [f (X)] : \cK (\hP, \P) \le \rho 
    \right\rbrace,
    \tag{\textsf{P}}
\end{align}
where $\rho \in [0, \infty)$, $\cP (\cX)$ is the set of all probability distributions on a data space $\cX$, $f: \cX \to \R$ is a loss function, $X$ is a random variable on $\cX$ having a nominal distribution $\hP$, and  {$\cK$ denotes the Kantorovich transport cost}, defined as 
\begin{align}
    \label{eqn:kantorovich dist}
    \cK (\hP, \P) = \inf _{\gamma \in \Gamma (\hP, \P)} \E _{(\hX, X) \sim \gamma} \big[ c (\hX, X) \big],
\end{align}
where $\Gamma (\hP, \P)$ denotes the set of all probability distributions on $\cX \times \cX$ with marginals $\hP$ and $\P$, and $c: \cX \times \cX \to [0, \infty]$ is a transport cost function.
 {Note that $\cK$ is a distance on $\cP (\cX)$ when $c$ is a metric.}
The function $\RP$ represents the robust loss hedging against deviations of data within $\rho$-neighborhood of the nominal distribution $\hP$.
When $c = d ^p$, where $d$ is a metric on $\cX$ and $p \in [1, \infty)$, problem \eqref{problem} is the inner worst-case problem in $p$-Wasserstein distributionally robust optimization, which has raised much interest recently; see \cite{kuhn2019wasserstein,blanchet2021statistical} for tutorials.

A central question of interest is developing the dual problem for \eqref{problem}.
In this paper, we present a novel proof that yields results held in more general settings than those found in the literature.
Assuming the interchangeability principle --- a condition inherent in existing results --- our proof is significantly shorter. The key idea is to view the robust loss as a function of the radius of the uncertainty set and then apply the Legendre transform to it twice. The concavity of the robust loss enables us to establish strong duality directly through the Legendre transformation. Compared to existing duality proofs that relied on convexity duality in conic programming \citep{esfahani2018data,zhao2018data,zhen2021mathematical} or vector spaces \citep{blanchet2019quantifying}, our proof only uses one-dimensional convex analysis. For detailed comparisons, we refer to Table \ref{tab:compare} and the discussion following Theorem \ref{thm}. Furthermore, we explore the interchangeability principle in depth and provide an equivalent condition, which has a strong connection to the measurable projection theorem and the measurable selection theorem. As a result, establishing the duality boils down to verifying the measurability conditions, which have been studied in broader settings.
To showcase its potential applications, we develop duality results for distributionally robust Markov decision processes, distributionally robust multistage stochastic programming, infinity-Wasserstein distributionally robust optimization, risk-averse optimization, and globalized distributionally robust counterpart.

\begin{table}
    \centering
    \caption{Comparison with existing duality results}\label{tab:compare}
    \small\setlength{\tabcolsep}{0pt}
    \resizebox{\textwidth}{!}{
    \begin{tabular}{c|cccccc}
        \toprule
        & \footnotesize\citet{esfahani2018data} &  \footnotesize\citet{zhao2018data}
        & \footnotesize\citet{blanchet2019quantifying} 
        
        & \footnotesize\citet{gao2022distributionally} 
        &\footnotesize\citet{zhen2021mathematical}
        & This paper
        \\\midrule
        $c$ & norm & continuous & lower semi-continuous  & power of metric & convex & arbitrary \\
        $f$ & piecewise concave & bounded & upper semi-continuous  & arbitrary & piecewise concave & arbitrary \\
        $\hP$ & empirical & empirical  & Borel   & Borel & empirical &
        \multirow{2}{*}{\hyperref[IP]{\shortstack{\textnormal\textsc(Interchangeability\\ Principle)}}} \\
        $\cX$   & \hspace{1pt} convex subset of $\R^d$\hspace{2pt} & convex compact    & Polish &\hspace{1pt}  Polish & convex subset of $\R^d$\hspace{2pt}  \\ 
        Proof\hspace{2pt} & conic duality & conic duality  & approximation argument  & constructive & Slater condition & Legendre transform\\
        \bottomrule
    \end{tabular}}
\end{table}

The remainder of this paper is structured as follows.
In Section \ref{sec:main}, we present our main proof. 
In Section \ref{sec:interchangeability}, we provide a verification result of the interchangeability principle.
We provide several examples in Section \ref{eg:generality} and extend our results to other distributional robust problems in Section \ref{sec:extensions}.
Finally, we conclude the paper in Section \ref{sec:conclusion}. 

\section{Model Formulation and Main Result}\label{sec:main}
In this section, we state and prove our main duality results under general assumptions. Proofs of auxiliary lemmas are provided in Appendix \ref{app:aux}.

\subsection{Main Assumptions}
\textbf{Notations.} We denote by $\Rbar := \R \cup \{ \pm \infty \}$ the extended reals and adopt the convention that $0 \cdot \infty = \infty$. 
 {For a probability space $(\cX,\sF, \hP)$}, we say an extended real-valued function on $(\cX,\sF)$ is measurable if it is $(\sF, \mathscr B (\Rbar))$-measurable, where $\mathscr B (\Rbar)$ is the Borel $\sigma$-algebra on $\Rbar$, and we say an extended real-valued function on $\cX$ is $\hP$-measurable if it is measurable with respect to the completion of $\sF$ under the measure $\hP$ \citep[Definition 1.11]{ambrosio2000functions}. 
We allow the expectation to take values in $\Rbar$ --- recall that the integration of a measurable function under a measure is well-defined whenever the positive part or the negative part of the integrand has a finite integral. 
When $\cX$ is a metric space, we denote by $d: \cX \times \cX \to [0, \infty)$ its metric.
Let $\cPs (\cX)$ denote the set of probability measures $\P$ on $(\cX, \sF)$ satisfying $\cK (\hP, \P) < \infty$.
 {For a function $h: \R \to \R \cup \{+\infty\}$, we denote by $h ^*: \R \to \R \cup \{+\infty\}$ its Legendre transform $h ^* (\lambda) := \sup _{\rho \in \R} \{ \lambda \rho - h (\rho) \} $. If $h: \R \to \Rbar$ attains $-\infty$ somewhere, then $h ^* \equiv +\infty$.}
 {
In addition to the problem \eqref{problem}, we also study its soft-penalty counterpart
\begin{align}\label{problem:soft}
    \sup _{\P \in \cPs} \left\lbrace
        \E _{X \sim \P} [ f (X) ] - \lambda \cK (\hP, \P)
    \right\rbrace,
    \tag{\textsf{P-soft}}
\end{align}
where $\lambda \in [0, \infty)$.}

We assume the following situation.
  
\begin{assumption}\label{assum:setup}
     {Let $(\cX, \sF, \hP)$ be a probability space}, $f: \cX \to \R$ be a measurable function with $\E _{\hP} [f] > -\infty$, and 
    $c: \cX \times \cX \to [0, \infty]$ be a measurable transport cost function satisfying $ {c(x, x) = 0}$ for all $x\in\cX$. 
\end{assumption}

The following lemma shows some useful properties of the worst-case loss $\RP (\cdot)$ defined in \eqref{problem}.

\begin{lemma}
    \label{lem:concave}
    Assume Assumption \ref{assum:setup} holds. Then
    $\RP (\cdot)$ is lower bounded by $\E_{\hP}[f]$, monotonically increasing, and concave on $[0,\infty)$. 
\end{lemma}

As we will see, the interchangeability principle discussed next is essential for strong duality, being both a necessary and sufficient condition. This principle is often mentioned in stochastic programming literature, as seen in references like \cite[Section 9.3.4]{shapiro2021lectures} and \cite{shapiro2017interchangeability}. It facilitates the swapping of expectation and minimization operators in our analysis.

\begin{IP}[Interchangability Principle \reflabel{IP}{\upshape (\textsf{IP})}]
    We say an $(\sF \otimes \sF)$-measurable function $\phi: \cX \times \cX \to \Rninfty$ satisfies the \emph{interchangeability principle} if the function $\hx\mapsto \sup _{x \in \cX} \phi (\hx, x)$ is $\hP$-measurable and it holds that 
    \begin{align*}
        \E _{\hX \sim \hP} \left[ 
            \sup _{x \in \cX} \phi (\hX, x) 
        \right]
        = 
        \sup _{\gamma \in \Gamma _{\hP}} \;
            \E _{(\hX, X) \sim \gamma} [ \phi (\hX, X) ],    
    \end{align*}
    where $\Gamma _{\hP}$ is the set of probability distributions on $(\cX \times \cX, \sF \otimes \sF)$ with first marginal $\hP$.
\end{IP}

\subsection{Main Result and its Proof}

Using Lemma \ref{lem:concave}, we derive the dual of \eqref{problem} as follows.

\begin{theorem}\label{thm}
    Assume Assumption \ref{assum:setup} holds. 
    Let $\lambda,\rho>0$.
    Then \eqref{problem:soft} is equivalent to
    \[
        (-\cL) ^* (-\lambda) = \sup _{\gamma \in \Gamma _{\hP}} 
        \E _{(\hX, X) \sim \gamma} \left[
            f (X) - \lambda c (\hX, X)
        \right],
    \]
    and 
    \eqref{problem} is equivalent to 
    \[
        \cL (\rho) = \min _{\lambda \ge 0} \left\lbrace
        \lambda \rho + \sup _{\gamma \in \Gamma _{\hP}}  \left\lbrace
            \E _{(\hX, X) \sim \gamma} \left[ 
                f (X) - \lambda c (\hX, X)
            \right]
            \right\rbrace
        \right\rbrace.
    \]
    In addition, for $\lambda > 0$, if and only if the function $\phi _{\lambda} (\hx, x) := f (x) - \lambda c (\hx, x)$ satisfies \ref{IP}, it holds that
    \begin{align*}
        (-\cL) ^* (-\lambda) = \E _{\hX \sim \hP} \left[ 
            \sup _{x \in \cX} \left\lbrace 
                f (x) - \lambda c (\hX, x) 
            \right\rbrace 
        \right],
    \end{align*}
    and if and only if $\phi_{\lambda}$ satisfies \ref{IP} for every $\lambda > 0$, it holds that
    \begin{align}
        \label{eqn:duality}
        \RP (\rho) = \min _{\lambda \ge 0} \left\lbrace
            \lambda \rho + \E _{\hX \sim \hP} \left[ 
                \sup _{x \in \cX} \left\lbrace 
                    f (x) - \lambda c (\hX, x) 
                \right\rbrace 
            \right]
        \right\rbrace, \quad \forall \rho > 0.
    \tag{\textsf{D}}
    \end{align}
\end{theorem}

\begin{remark}[$p$-Wasserstein distance]
    Recall that the $p$-Wasserstein distance $\cW _p (\hP, \P)$, $p\in[1 ,\infty)$, is defined by 
    \begin{align*}
        \cW _p (\hP, \P) = \inf _{\gamma \in \Gamma (\hP, \P)} \|d\| _{L ^p (\cX \times \cX; \gamma)} =
            \inf _{\gamma \in \Gamma (\hP, \P)} \E _{(\hX, X) \sim \gamma} [d (\hX, X) ^p] ^\frac1p.
    \end{align*}
    By setting $c (\hx, x) = d (\hx, x) ^p$ we have $\cK (\hP, \P) = \cW _p ^p (\hP, \P)$. Thereby \eqref{eqn:duality} corresponds to the dual formulation of the $p$-Wasserstein DRO
    \begin{align*}
        \cL (\rho ^p) = \min _{\lambda \ge 0} \left\lbrace
        \lambda \rho^p + 
            \E _{\hX\sim\hP} \left[ \sup _{x\in\cX} \left\{
                f (X) - \lambda d (\hX, x) ^p
                \right\}
            \right]
        \right\rbrace.
    \end{align*}
    We will handle the case $p = \infty$ separately in Section \ref{sec:maximum-cost}.
\QDEF
\end{remark}

\begin{remark}[Necessity of \ref{IP}]\label{rmk:M}
    The second part of Theorem \ref{thm} discusses the necessity of the interchangeability principle \ref{IP}.
    As will be elaborated on in the next section, it ensures the measurability of the supremum function in \eqref{eqn:duality} and the existence of approximately worst-case distributions. To our knowledge, \ref{IP} is weaker than the assumptions in all existing results that enable the expression \eqref{eqn:duality}.
\QDEF    
\end{remark}

\begin{remark}[Continuity at $\rho=0$]
    In general, \eqref{eqn:duality} does not hold at $\rho = 0$.
    Indeed, the right-hand side of \eqref{eqn:duality} is continuous in $\rho \in [0, \infty)$, but $\cL (\rho)$ may be not right-continuous at $0$. 
    For instance, if $\cX = \R$, $c (\hx, x) = |\hx - x|$, $f (x) = \mathbf1{\{x \neq 0\}}$, and $\hP = \Dirac _0$, the Dirac measure at $0$, then $\cL (\rho) = 1$ for any $\rho > 0$ and $\cL (0) = 0$. A sufficient condition ensuring the right-continuity of $\cL (\rho)$ at $0$ is the following: there exists a continuous concave function $\varphi: [0, \infty) \to [0, \infty)$ with $\varphi (0) = 0$ such that
    $
        f (x) - f (\hx) \le \varphi \circ c (\hx, x)
    $
    for all $x \in \cX$ and $\hP$-a.e. $\hx \in \cX$ with $c (\hx, x) < \infty$.
    Indeed, under this condition, for any $\epsilon > 0$ and $ {\P \in \cP (\cX)}$ with $\cK (\hP, \P) \le \epsilon$, there exists a $\gamma \in \Gamma _{\hP}$ such that $\E _{\gamma} [c (\hX, X)] \le 2 \epsilon$, hence $\E _{\gamma} [f (X) - f (\hX)] \le \E _\gamma [\varphi \circ c] \le \varphi (\E _{\gamma} [c])$ by Jensen's inequality. Therefore,
    $
        \cL (\epsilon) \le \E _{\hP} [f] + \varphi (\E _{\gamma} [c])  {\le} \cL (0) + \varphi (\E _{\gamma} [c])
    $, which converges to $\cL (0)$ as $\epsilon \to 0$.
    When $c = d ^p$, where $d$ is a metric on $\cX$, this condition is related to the growth condition imposed in \cite{gao2022distributionally} and the upper semi-continuity of $f$.
\QDEF
\end{remark}

\proof{Proof of Theorem \ref{thm}.}
      {Fix $\lambda > 0$. By definition, $\cL (\rho) = -\infty$ for $\rho < 0$.} Taking the Legendre transform of $-\RP (\cdot)$ gives that 
    \begin{align*}
        (-\cL) ^* (-\lambda) 
        & = \sup _{\rho \ge 0} \lbrace (-\lambda) \rho - (-\RP (\rho)) \rbrace \\
        &= \sup _{\rho \ge 0} \lbrace \RP (\rho) - \lambda \rho \rbrace \\
        &= \sup _{\rho \ge 0} \sup _{\P \in \cP(\cX)} \left\lbrace
            \E _{X \sim \P} [ f (X) ] - \lambda \rho:\; 
            \cK (\hP, \P) \le \rho
        \right\rbrace \\
        &= \sup _{\P \in \cP(\cX)} \sup _{\rho \ge 0} \left\lbrace
            \E _{X \sim \P} [ f (X) ] - \lambda \rho:\;  
            \cK (\hP, \P) \le \rho
        \right\rbrace \\
        &= \sup _{\P \in \cPs} \left\lbrace
            \E _{X \sim \P} [ f (X) ] - \lambda \cK (\hP, \P)
        \right\rbrace,
    \end{align*}
    which gives \eqref{problem:soft}.
    Using the definition of the  {Kantorovich transport cost} \eqref{eqn:kantorovich dist}, it follows that 
    \begin{align*}
        (-\cL) ^* (-\lambda)
        & = \sup _{\P \in \cPs} \left\lbrace
            \E _{X \sim \P} [ f (X) ] - \lambda \cK (\hP, \P)
        \right\rbrace \\
        &= 
        \sup _{\P \in \cPs} \left\lbrace
            \E _{X \sim \P} [ f (X) ] - \lambda \inf _{\gamma \in \Gamma (\hP, \P)} \E _{(\hX, X) \sim \gamma} [c (\hX, X)]
        \right\rbrace 
        \\
        &= 
        \sup _{\P \in \cPs, \gamma \in \Gamma (\hP, \P)} \left\lbrace
            \E _{X \sim \P} [ f (X) ] - \lambda \E _{(\hX, X) \sim \gamma} [c (\hX, X)]: 
            \E _{\gamma} [c] < \infty
        \right\rbrace 
        \\
        &=
        \sup _{\gamma \in \Gamma _{\hP}} \;
            \E _{(\hX, X) \sim \gamma} \left[
                f (X) - \lambda c (\hX, X)
            \right].
    \end{align*} 
     {Observe that $(-\cL) ^* (-\lambda) \ge \sup _{\rho \ge 0} \{(-\lambda) \rho + \cL (0)\} = +\infty$ for $\lambda < 0$.}
    From Lemma \ref{lem:concave} we have that $\RP (\cdot)$ is bounded from below, increasing and concave in $[0, \infty)$, so either $\cL (\rho) < +\infty$ for every $\rho \ge 0$ or $\cL (\rho) = +\infty$ for every $\rho > 0$.
    In the former case, by the involution property of Legendre transform (see, e.g., \cite[Theorem 12.2]{rockafellar1970}),  {$(-\cL) ^{**}$ equals to the lower semi-continuous convex envelope of $-\cL$. Since $\cL (\cdot)$ is concave in $[0, \infty)$, $-\cL (\rho) = -\cL ^{**} (\rho)$ on the interior of the set $\{\rho\in\R:-\cL(\rho) < +\infty\}$, which is $(0, \infty)$.}  
    Hence, for every $\rho > 0$,
    \begin{align*}
        \cL (\rho) = -(-\cL) ^{**} (\rho) &= -\max _{\lambda \ge 0} \left\lbrace
            (-\lambda) \rho - (-\cL) ^* (-\lambda)
        \right\rbrace \\
        &= 
        \min _{\lambda \ge 0} \left\lbrace
            \lambda \rho + (-\cL) ^* (-\lambda)
        \right\rbrace \\
        &= \min _{\lambda \ge 0} \left\lbrace
            \lambda \rho + \sup _{\gamma \in \Gamma _{\hP}}  \left\lbrace
                \E _{(\hX, X) \sim \gamma} \left[ 
                    f (X) - \lambda c (\hX, X)
                \right]
            \right\rbrace
        \right\rbrace.
    \end{align*}
     {Note that we switched from supremum to maximum because $(-\cL) ^* (-\lambda)$ is lower semi-continuous and bounded from below for $\lambda \ge 0$ by Lemma \ref{lem:convex}. Therefore, $(-\lambda) \rho - (-\cL) ^* (-\lambda)$ can be arbitrarily small as $\lambda \to +\infty$, hence the maximum is attainable. In the latter case, $(-\cL) ^* (-\lambda) = +\infty$ for any $\lambda \in \R$, so the above is also true.
    This proves the first part of the theorem.
    }  

     {
    For the second part, note that for $\lambda > 0$, $\phi _\lambda = f - \lambda c$ satisfies \ref{IP} means
    \begin{align*}
        (-\cL) ^* (-\lambda) = \sup _{\gamma \in \Gamma _{\hP}} \left\lbrace
            \E _{(\hX, X) \sim \gamma} [ f (X) - \lambda c (\hX, X) ]
        \right\rbrace = \E _{\hX \sim \hP} \left[ 
            \sup _{x \in \cX} \left\lbrace 
                f (x) - \lambda c (\hX, x) 
            \right\rbrace 
        \right] =: \cG (\lambda),
    \end{align*}
    where $\cG (\lambda)$ is defined as above for $\lambda \ge 0$ and $\cG (\lambda) := (-\cL) ^* (-\lambda) = +\infty$ for $\lambda < 0$.
    By Lemma \ref{lem:convex} in Appendix \ref{app:aux}, $(-\cL) ^* (-\lambda)$ and $\cG (\lambda)$ are lower bounded by $\E _{\hP} [f]$, monotonically decreasing, convex, and lower semi-continuous on $[0, \infty)$, implying the right continuity at $\lambda = 0$. Hence, the following equivalence holds:
    \begin{align*}
        & \phi _\lambda \text{ satisfies \ref{IP} for all } \lambda \in (0, \infty) \\
        \iff & (-\cL) ^* (-\lambda) = \cG (\lambda) \text{ for all } \lambda \in (0, \infty) \\
        \iff & (-\cL) ^* (-\lambda) = \cG (\lambda) \text{ for all } \lambda \in [0, \infty) \\
        \iff & (-\cL) ^* (\lambda) = \cG (-\lambda) \text{ for all } \lambda \in \R
    \end{align*}
    Below we separately consider the cases $\cL (\rho) < +\infty$ and $\cL (\rho) \equiv +\infty$ for $\rho \ge 0$.
    If $\cL (\rho) < +\infty$ for every $\rho \ge 0$, then $(-\cL) ^* \not\equiv +\infty$ by Lemma \ref{lem:convex-2}, so by the involution property of Legendre transform, we have
    \begin{align*}
        & (-\cL) ^* (\lambda) = \cG (-\lambda) \text{ for all } \lambda \in \R \\
        \iff & (-\cL) ^{**} (\rho) = (\cG (-\cdot)) ^* (\rho) = \cG ^* (-\rho)  \text{ for all } \rho \in \R \\
        \iff & (-\cL) ^{**} (\rho) = \cG ^* (-\rho) \text{ for all } \rho \in [0, \infty) \\
        \iff & (-\cL) ^{**} (\rho) = \cG ^* (-\rho) \text{ for all } \rho \in (0, \infty) \\
        \iff & \cL (\rho) = \cG ^* (-\rho) = \inf _{\lambda \in \R} \left\{ \lambda \rho + \cG (\lambda) \right\} = \min _{\lambda \ge 0} \left\{ \lambda \rho + \cG (\lambda) \right\} \text{ for all } \rho \in (0, \infty)
    \end{align*}
    Here we have used the proporties that $(-\cL) ^{**}$ and $\cG ^* (-\cdot)$ are both right continuous at $0$ and are both $+\infty$ on $(-\infty, 0)$, thanks to Lemma \ref{lem:convex-2}. In the last step, the minimum is attainable because $\lambda \rho + \cG (\lambda)$ is bounded from below and lower semi-continuous, and becomes arbitrarily large as $\lambda \to \infty$.
    If $\cL (\rho) \equiv +\infty$ for $\rho > 0$, then $(-\cL) ^* (-\lambda) = +\infty$ for all $\lambda \in \R$, so 
    \begin{align*}
        & (-\cL) ^* (\lambda) = \cG (-\lambda) \text{ for all } \lambda \in \R \\
        \iff & \cG (\lambda) = +\infty \text{ for all } \lambda \in \R
        \\
        \iff & \cG (\lambda) = +\infty \text{ for all } \lambda \ge 0
        \\
        \iff &
        \min _{\lambda \ge 0} \left\{ \lambda \rho + \cG (\lambda) \right\} = +\infty \text{ for all } \rho \in (0, \infty)
        \\
        \iff &
        \min _{\lambda \ge 0} \left\{ \lambda \rho + \cG (\lambda) \right\} = \cL (\rho) \text{ for all } \rho \in (0, \infty)
    \end{align*}
    In conclusion, $\phi _\lambda$ satisfies \ref{IP} for all $\lambda \in (0, \infty)$ if and only if $\cL (\rho) = \min _{\lambda \ge 0} \left\{ \lambda \rho + \cG (\lambda) \right\}$ for all $\rho \in (0, \infty)$.}
\QED
\endproof

Let us compare our proof technique with existing duality results in the literature.
Proofs in \citep{esfahani2018data,blanchet2019quantifying,zhao2018data,sinha2018certifying,zhen2021mathematical} rely on advanced convex duality theory.
More specifically,
\citet{esfahani2018data,zhao2018data,zhen2021mathematical} exploit conic duality \citep{shapiro2001duality} for the problem of moments. This approach requires the nominal distribution $\hP$ to be finitely supported and the space $\cX$ to be convex, along with additional assumptions on the transport cost $c$ and the loss function $f$.
\citet{blanchet2019quantifying} uses an approximation argument that represents the Polish space $\cX$ as an increasing sequence of compact subsets. This enables duality for any Borel distribution $\hP$, based on the Fenchel conjugate on vector spaces \citep{luenberger1997optimization}, under semi-continuity assumptions on the transport cost function $c$ and loss function $f$.
Using the same infinite-dimensional convex duality, \citet[Theorem 5]{sinha2018certifying} streamlines the analysis by assuming the function $(\hX,X) \mapsto \lambda c(\hX,X) - f(X)$ is a normal integrand \citep{rockafellar2009variational}.
Compared with these non-constructive duality proofs, our (non-constructive) proof employs only the Legendre transform, namely, the convex duality for univariate real-valued functions.
The constructive proof developed by \citet{gao2022distributionally} achieves a similar level of generality as \cite{blanchet2019quantifying}, but without relying on convex duality theory. They construct an approximately worst-case distribution using the first-order optimality condition of the weak dual problem.
Although both their approach and ours avoid using advanced minimax theorems, our analysis is notably shorter.

\section{Discussion on the Interchangeability Principle}\label{sec:interchangeability}

In this section, we first show that the interchangeability principle \ref{IP} in Section \ref{sec:main} is intricately linked to the measurable projection and measurable selection conditions,  {and then prove \ref{IP} holds in the Wasserstein DRO setting}.

 {
Note that in \eqref{eqn:duality}, the function $f$ belongs to a family of loss functions, and the transport cost function $c$ is chosen contingent on the specific application.  Consequently, from a pragmatic standpoint, we aim for \ref{IP} to be applicable to functions within the family 
\begin{align*}
    \left\lbrace 
        \phi: \cX \times \cX \to \R \cup \{-\infty\} \,\bigg\vert\, 
        \begin{array}{ll}
            \phi (\hx, x) = f (x) - \lambda c (\hx, x),\ f: \cX \to \R \text{ measurable, } \lambda \ge 0 \\
            c: \cX \times \cX \to [0, \infty] \text{ measurable, }c (x, x) = 0 \text{ for all } x \in \cX
        \end{array}
    \right\rbrace
\end{align*}
Observe that the non-negativity and reflexivity of $c$ imply the following \emph{diagonally dominant} property of functions in this family (see Lemma \ref{lem:f-lc} in Appendix \ref{app:aux}).
}

 {
\begin{definition}
    \label{def:diagonally-dominant}
    We say a $(\sF \otimes \sF)$-measurable function $\phi: \cX \times \cX \to \Rninfty$ is \emph{diagonally dominant} if $\phi (\hx, x) \le \phi (x, x)$ for every $\hx, x \in \cX$. 
    We say a $(\sF \otimes \sF)$-measurable set $A \subset \cX \times \cX$ is \emph{diagonally dominant} if for every $(\hx, x) \in A$, it holds that $(x, x) \in A$. 
    We say a set function $E: \cX \to \sF \setminus \lbrace\varnothing\rbrace$ with $(\sF \otimes \sF)$-measurable graph is \emph{diagonally dominant} if $x \in E (\hx)$ implies $x \in E (x)$.
    \QDEF
\end{definition}
}

With this definition, we show that \ref{IP} is equivalent to the \emph{measurable projection} and the weak \emph{measurable selection} conditions.
We denote by $\Fcomp$ the completion of $\sF$ under $\hP$ \citep[Page 13]{kallenberg1997foundations}.

\begin{proposition}
    \label{prop:opt-selection}
    \ref{IP} holds for all $(\sF \otimes \sF)$-measurable diagonally dominant functions $\phi: \cX \times \cX \to \Rninfty$ if and only if
    $(\cX, \sF, \hP)$ satisfies the following two conditions:
    \begin{enumerate}[leftmargin=*]
        \item [\reflabel{ass:Proj}{\upshape (\textsf{Proj})}]
        {\upshape[Measurable Projection]} 
        For any  {diagonally dominant} set $A \in \sF \otimes \sF$, 
        \begin{align*}
            \Proj _{\hx} (A) := \left\lbrace \hx \in \cX: (\hx, x) \in A \text{ for some } x \in \cX \right\rbrace \in \Fcomp.
        \end{align*}

        \item [\reflabel{ass:Sel*}{\upshape (\textsf{Sel*})}]
        {\upshape[Weak Measurable Selection]} 
        For any  {diagonally dominant} set-valued function $E: \cX \to \sF \setminus \lbrace\varnothing\rbrace$ with a measurable graph
        \begin{align*}
            \Graph (E) := \left\lbrace (\hx, x) \in \cX \times \cX: x \in E(\hx) \right\rbrace \in \sF \otimes \sF,
        \end{align*}
        there exists a probability measure $\gamma \in \Gamma _{\hP}$, such that $\supp \gamma \subset \Graph (E)$.
    \end{enumerate}
\end{proposition}

\begin{remark}
    Measurable projection and measurable selection are often seen in the literature on stochastic control (e.g., \cite{rockafellar2009variational,bertsekas1996stochastic}) and stochastic programming (e.g., \cite{shapiro2017interchangeability}).
    \ref{ass:Proj} states that the projection operator $\Proj _{\hx}: (\hx, x) \mapsto \hx$ maps measurable sets in $\sF\otimes\sF$ to $\hP$-measurable sets. 
    We refer \ref{ass:Sel*} to as the weak measurable selection, because it is weaker than the measurable selection in the literature \cite{rockafellar2009variational,shapiro2021lectures}.
    The measurable selection therein involves a \textsl{deterministic} selection which, in our context, reads as
    \begin{enumerate}[leftmargin=*]
        \item [\reflabel{ass:Sel}{\upshape (\textsf{Sel})}]
        [Measurable Selection] \emph{
        For any set-valued function $E: \cX \to \sF \setminus \lbrace\varnothing\rbrace$ with a measurable graph,
        there exists an $(\Fcomp, \sF)$-measurable map $T: \cX \to \cX$ such that $T (\hx) \in E (\hx)$, $\forall \hx \in \cX$.}
    \end{enumerate}
    \citet[Theorem 14.60]{rockafellar2009variational} (see also \cite[Section 9.3.4]{shapiro2021lectures}) shows that
    measurable projection and measurable selection together imply the following
    \[
        \E _{\hX \sim \hP} \left[
            \sup _{x \in \cX} \phi (\hX, x)
        \right] = \sup _{T \in \cT} \E _{\hX \sim \hP} \left[
            \phi (\hX, T (\hX))
        \right],
    \]
    where $\cT$ denotes the set of $(\sF _{\hP}, \sF)$-measurable maps. Note that $(\Id \otimes T) _\# \hP \in \Gamma _{\hP}$, so the above is stronger than \ref{IP}.
    Comparatively, the weak measurable selection condition \ref{ass:Sel*} allows a \textsl{random} selection, represented by the conditional distribution of $\gamma _{x | \hx}$ supported on $E (\hx)$. This indicates that \ref{ass:Sel*} is weaker than \ref{ass:Sel}.
    By Proposition \ref{prop:opt-selection}, \ref{IP} is equivalent to \ref{ass:Sel*} and \ref{ass:Proj}.
\QDEF
\end{remark}

The next result shows that \ref{IP} holds when the transport cost function corresponds to the $p$-Wasserstein DRO, even if $f$ is not $\sF$-measurable but merely $\hP$-measurable. This will be used in Example \ref{eg:MDP}.

\begin{proposition}
    \label{prop:c-continuous}
    Let $(\cX, d)$ be a metric space equipped with Borel $\sigma$-algebra $\sF$, $\hP$ be a tight measure, $f$ be $\hP$-measurable with $\E _{\hP} [f] > -\infty$. Let $p \in [1,\infty)$, $\lambda \ge 0$. Then the function $\phi (\hx, x) = f (x) - \lambda d (\hx, x) 
    ^p$ is $(\sF \otimes \sF _{\hP})$-measurable and satisfies \ref{IP}.
\end{proposition}

\section{Examples}\label{eg:generality}
In this section, we offer several examples that demonstrate how our findings not only align with existing research but also are useful for important applications in the area of distributionally robust sequential decision-making.

The following two examples illustrate that existing results in the literature are based on assumptions that are strictly stronger than \ref{ass:Proj} and \ref{ass:Sel*}. Consequently, our results strictly generalize existing findings in the literature.

\begin{example}[Empirical distribution]\label{eg:empirical}
    If $(\cX, \Fcomp)$ is a discrete measurable space, that is, $\Fcomp = 2 ^\cX$, then \ref{ass:Proj}  
    and \ref{ass:Sel} 
    always holds, because every subset of $\cX$ is $\Fcomp$-measurable,
    and every map $T: \cX \to \cX$ is $\hP$-measurable.
    For instance, when $\cX$ is a measurable space equipped with a Borel $\sigma$-algebra and $\hP$ is finitely supported, then $\sF _{\hP} = 2 ^{\cX}$ is the collection of all subsets of $\cX$, which is the discrete $\sigma$-algebra on $\cX$.
    Thus our result covers the results in \cite{esfahani2018data,zhao2018data,zhen2021mathematical}, which studies the case where $\cX$ is a convex subset of $\R^d$. 
    \QEG
\end{example}

\begin{example}[Polish space]\label{eg:polish}
    If $\cX$ is a Polish (complete and separable metric) space and $\sF$ is its Borel $\sigma$-algebra, then \ref{ass:Proj} holds due to \cite[Theorem 8.3.2]{aubin2009set}, and \ref{ass:Sel} holds due to \cite[Theorem 18.26]{aliprantis2006infinite}.
    Thereby our result covers the results in \cite{blanchet2019quantifying,gao2022distributionally}.
    \QEG
\end{example}

Example \ref{eg:polish} can be generalized as follows.

\begin{example}[Suslin space]\label{eg:suslin}
    A Hausdorff topological space is Suslin (also known as analytic) if it is the continuous image of a Borel set in a Polish space.
    If $\cX$ is a Suslin space and $\sF$ is its Borel $\sigma$-algebra, then \ref{ass:Proj} holds due to \cite[Theorem III.23]{castaing1977measurable}, 
    and \ref{ass:Sel} holds due to \cite[Theorem III.22]{castaing1977measurable}.
    \QEG
\end{example}

The following examples showcase the broad applicability of our results.

\newcommand{\Kappa}{\mathrm{K}}
\newcommand{\hs}{\widehat{s}}
\begin{example}[Distributionally robust Markov decision process]\label{eg:MDP}
    Consider a finite-horizon Markov decision process.
    The standard working horse \citep{bertsekas1996stochastic,shreve1979universally} involves the Borel space (i.e., a Borel subset of a complete and separable metric space).
     {To avoid measurability issues, existing literature often assumes a finite or countable state space. 
    Nonetheless, for general Borel state space, we can still verify \ref{IP} by Proposition \ref{prop:c-continuous}.}
    Let the state space $\cS$ and the action space $\cA$ be non-empty Borel spaces. 
    Let $\{\cA (s)\} _{s \in \cS} \subset \cA$ be a family of non-empty feasible action sets such that the corresponding constraint set
    $\Kappa = \{(s, a): s \in \cS, a \in \cA (s)\}$ is an analytic subset of $\cS\times\cA$.
    Let the one-stage cost $g_t$ be a lower semi-analytic function on $\Kappa$ bounded from below.
    Suppose we have a nominal transition kernel $\{\hP (\cdot \mid s, a)\} _{(s, a) \in \Kappa}$ that is a Borel measurable stochastic kernel on $\cS$ given $(s, a)$.
    Consider the following uncertainty sets defined for every state-action pair $(s, a)$
    \begin{align*}
        \frakM (s, a) = \left\{\P \in \cP (\cS): \cK (\hP (\cdot \mid s, a), \P) \le \rho (s, a) \right\}, \qquad (s, a) \in \Kappa,
    \end{align*}
    where the positive radius function $\rho$ is lower semi-analytic on $\Kappa$, and the transport cost $c$ associated with $\cK$  {is $d ^p$, $p \in[1,\infty)$ where $d$ is the metric on $\cS$}.
    The distributionally robust counterpart of the value iteration \citep{yang2017convex,wang2020reliable} is given by $V _{T + 1} \equiv 0$, and for $t = 1, \ldots, T$,
    \begin{align*}
        V _t (s) = \inf _{a \in \cA (s)} \left\{
            g _t (s, a) + \sup _{\P \in \frakM (s, a)} \E _{s' \sim \Pr} [V _{t + 1} (s')] 
        \right\}.
    \end{align*}
    Then by induction, we can prove the duality 
    \begin{align}
        V _t (s) = \inf _{\substack{a\in\cA(s) \\ \lambda \ge 0}} \left\{ 
            g _t (s, a) + \lambda \rho (s, a) + \E _{\hs' \sim \hP (\cdot \mid s, a)} \left[
                \sup _{s' \in \cS} \big\{ 
                    V _{t + 1} (s') - \lambda c (\hs', s')
                \big\} 
            \right] 
        \right\}.
    \end{align}
     {In the inductive step, $V _{t + 1}$ is $\hP$-measurable and lower semi-analytic. The strong duality holds thanks to Proposition 2. Since $\{\hP (\cdot \mid s, a)\} _{(s, a) \in \Kappa}$ is a Borel measurable kernel, the robust loss is also lower semi-analytic. Taking the infimum yields a lower semi-analytic value function $V _t$, which completes the induction. We refer to Appendix \ref{app:generality} for details.}
    \QEG
\end{example}

\begin{example}[Data-driven robust multi-stage stochastic programming]\label{eg:SP}
    Consider a multi-stage stochastic programming problem \citep{shapiro2021lectures}.
    Let $(x _1, \ldots, x _T)$ be a $T$-stage random data process, which is assumed to be stagewise independent, namely, $\{x _t\} _{t = 1} ^T$ are mutually independent. At each stage, after the current-stage uncertainty $x _t$ is realized, the decision maker seeks a non-anticipative decision $u _t$ from a feasible set $\cU _t (u _{t - 1}, x _t) \subset \R ^{d _t}$ dependent on the previous-stage decision $u _{t - 1}$ and the current-stage uncertainty $x _t$, where $\cU _t$ is a measurable multifunction (i.e., set-valued function).
    The (random) cost of taking a decision $u _t$ at stage $t$ is $f _t (u _t, x _t)$ and the goal is to minimize the cumulative cost.
    Suppose the uncertainty $x _t$ has an unknown distribution on a Suslin space $\cX _t$ equipped with the Borel $\sigma$-algebra, and one formulates an uncertainty set based on a nominal Borel distribution $\hP _t$ with a radius $\rho _t>0$:
    \begin{align*}
        \frakM _t = \left\{
            \P \in \cP (\cX _t): \cK (\hP _t, \P) \le \rho _t 
        \right\}, \qquad t = 1, \ldots, T.
    \end{align*}
    Following the convention in stochastic programming, we assume there is no uncertainty in the first stage, i.e., $\rho _1 = 0$ and $\hP _1$ is a Dirac measure $\hP _1 = \Dirac _{x _1}$, where $x _1 \in \cX _1$, and we set $u _0 = \varnothing$.
    It would be natural to consider the following distributionally robust counterpart of Bellman recursion: $Q _{T + 1} \equiv 0$, and for $t = 1, \ldots, T$,
    \begin{align*}
        Q _t (u _{t - 1}, x _{t}) = \inf _{u \in \cU _t (u _{t - 1}, x _t)} \left\{
            f _t (u, x _t) + \sup_{\P \in \frakM _{t + 1}} \E_{\P} \big[
                Q _{t + 1} (u, x _{t + 1})
            \big]
        \right\}.  
    \end{align*}
    Under this setting, \ref{ass:Proj} and \ref{ass:Sel} follow from Example \ref{eg:suslin}.    
    Assume that $f _t: \R ^{d _t} \times \cX _t \to \R$ is bounded from below and random lower semi-continuous, i.e. the epigraph multifunction $x _t \mapsto \operatorname{epi} f _t (\cdot, x _t)$ is measurable and closed valued \cite[Definition 9.47]{shapiro2021lectures}. Assume the multifunction $\cU _t: \R ^{d _{t - 1}} \times \cX _t \rightrightarrows \R ^{d _t}$ is measurable. Moreover, for each $x _t \in \cX _t$, assume $\cU _t (\cdot, x _t)$ is non-empty and has closed graph. Assume further that there exists a bounded set containing $\cU _t (u _{t - 1}, x _t)$ for all $u _{t - 1}$ and $x _t$.
    Then by induction (see Appendix \ref{app:generality} for details), we have the duality
    \begin{align*}
        Q _t (u _{t - 1}, x _t) = \inf _{\substack{u \in \cU _t (u _{t - 1}, x _t) \\ \lambda \ge 0}} \left\{
            f _t (u, x _t) + \lambda \rho _t + \E _{\hP _{t + 1}} \left[
                \sup _{x \in \cX _{t + 1}} \big\{
                    Q _{t + 1} (u, x) - \lambda c (\hx _{t + 1}, x)
                \big\} 
            \right] 
        \right\}.
        \tag*{\QEG} 
    \end{align*}
\end{example}

\begin{example}[Chance constraint]
    \label{eg:01}
    Let $(\cX, d)$ be a metric space and let $\sF$ be its Borel $\sigma$-algebra.
    Consider a distributionally robust chance constraint \cite{chen2023approximations,xie2021distributionally,yang2022wasserstein}
    \begin{align}\label{eqn:cc}
        \inf _{\P \in \cP (\cX)} \left\lbrace 
            \P (\cS): \cW _p (\hP, \P) \le \rho 
        \right\rbrace 
        \ge 1 - \beta,
    \end{align}
    which ensures that an event $\cS \in \sF$ happens with probability higher than $1 - \beta$, $\beta \in (0, 1)$, with respect to every distribution within the uncertainty set. 
    Denote by $\cSc$ the complement of the set $\cS$.
    In Appendix \ref{app:generality}, we show that if $\phi (\hx, x) = \mathbf1 _{\cSc} (x) - \lambda d (\hx, x) ^p$ satisfies \ref{IP} for every $\lambda > 0$, then using Theorem \ref{thm}, the distributionally robust chance constraint \eqref{eqn:cc} is equivalent to a worst-case conditional value-at-risk constraint 
    \begin{align*}
        \rho ^p \le -\beta \CVaRhP _{\beta} (-d (\hX, \cSc) ^p).
    \end{align*}
    where $\CVaRhP _\beta(\cdot)$ represents the conditional value-at-risk at risk level $\beta$.
    Note that the constraint is infeasible if $\rho \ge \E _{\hX \sim \hP} [d (\hX, \cSc) ^p] ^{\frac1p}$.
    Similar results have been obtained by \cite{chen2023approximations,yang2022wasserstein,xie2021distributionally} but with less generality.
    In \cite{chen2023approximations}, the ambient space $\cX$ is Euclidean, and $\hP$ is empirical. In \cite{yang2022wasserstein}, the ambient space is a subset of a Euclidean space and $\hP$ is a Borel probability measure. In \cite{xie2021distributionally}, $\cX$ is a totally bounded Polish space. Both \cite{chen2023approximations,xie2021distributionally} consider only 1-Wasserstein distance. In comparison, the result above requires $(\cX, d)$ only to be a metric space satisfying \ref{IP} and has no further requirement on the nominal probability distribution $\hP$ nor the set $\cS$ beyond measurability.
    \QEG
\end{example}

\section{Extensions}
\label{sec:extensions}

In this section, we extend our results to several other problems. The proofs are based on similar techniques that we developed in Section \ref{sec:main}.

\subsection{Infinity-Wasserstein DRO and Maximum Transport Cost}
\label{sec:maximum-cost}

Recall the $\infty$-Wasserstein distance
\[
    \cW _\infty (\hP, \P) = \inf _{\gamma \in \Gamma (\hP, \P)} \|d\| _{L ^\infty (\cX \times \cX; \gamma)}
    = \inf _{\gamma \in \Gamma (\hP, \P)} \gamma \text{-} \essssup _{\hx, x \in \cX} d (\hx, x),
\]
The result in Section \ref{sec:main} only covers the Wasserstein distance of a finite order.
To study the case $p = \infty$, we introduce the \emph{maximum transport cost} as 
\begin{align*}
    \cM (\hP, \P) := \inf _{\gamma \in \Gamma (\hP, \P)} \gamma \text{-} \essssup _{\hx, x \in \cX} c (\hx, x),
\end{align*}
where $c$ is a transport cost function satisfying Assumption \ref{assum:setup}.
We define the maximum transport cost robust loss by 
\begin{align*}
    \cLm (\rho) &:= \sup _{\P \in \cP (\cX)} \left\lbrace
        \E _{X \sim \P} [f (X)]: \cM (\hP, \P) \le \rho 
    \right\rbrace.
\end{align*}
Similar to the results in Section \ref{sec:main}, the soft-constrained counterpart can be viewed as the Legendre transform of negative hard-constrained robust loss
\[
  (-\cLm) ^* (-\lambda) = \sup _{\rho \ge 0} \left\lbrace (-\lambda) \rho - (-\cLm (\rho)) \right\rbrace = 
    \sup _{\P \in \cP (\cX)} \left\lbrace
        \E _{X \sim \P} [f (X)] - \lambda \cM (\hP, \P) 
    \right\rbrace.
\]
With this definition, we cover the $\infty$-Wasserstein DRO by setting $c = d$.
Below, we first establish a duality result when the constraint of the uncertainty set is a strict inequality.

\begin{proposition}
    \label{prop:infty-<}
    Let $f$ and $c$ satisfy Assumption \ref{assum:setup}.
    Define
    \[
      \cLm ^\circ (\rho) := \sup _{\P \in \cP (\cX)} \left\lbrace
            \E _{X \sim \P} [f (X)]: \cM (\hP, \P) < \rho 
        \right\rbrace.
    \]
    Suppose the function $\psi _\rho (\hx, x) := f (x) - \infty \mathbf1{\{c (\hx, x) \ge \rho\}}$ satisfies \ref{IP} for some $\rho > 0$. Then
    \begin{align*}
        \cLm ^\circ (\rho) 
        &= \E _{\hX \sim \hP} \left[\sup _x \left\lbrace
            f (x): c (\hX, x) < \rho
        \right\rbrace \right].
    \end{align*}
    Suppose the $\psi _\rho$ satisfies \ref{IP} for every $\rho > 0$. Then
    \begin{align*}
        (-\cLm) ^* (-\lambda) 
        = \sup _{\rho > 0} \left\lbrace 
        \E _{\hX \sim \hP} \left[\sup _x \left\lbrace
            f (x): c (\hX, x) < \rho 
        \right\rbrace - \lambda \rho
        \right]
    \right\rbrace.
    \end{align*}
\end{proposition}

\begin{remark}
    Unlike the results in Section \ref{sec:main}, we have to distinguish two cases: $\cLm ^\circ (\rho)$ which involves the strict inequality constraint, and $\cLm (\rho)$ which involves the non-strict inequality constraint. 
    Indeed, for problem \eqref{problem} in Section \ref{sec:main}, the value of $\cL (\rho)$ in \eqref{problem} would not be affected if we replace the non-strict inequality by the strict inequality thanks to the concavity, and thus continuity, of $\cL$ with respect to $\rho \in (0, \infty)$. However, for the problem with the maximum transport cost, neither $\cLm ^\circ$ nor $\cLm$ is necessarily continuous. 
    
    Without additional assumptions, the duality does not hold for the equality constraint, even when the cost function $c$ is also a metric. For instance, consider $\cX = [0, 1] \times [0, 1]$, $\hP$ is a uniform distribution on $\{0\} \times [0, 1]$, and $f (x _1, x _2) = \mathbf1{\{x _1 = 1\}}$. We introduce the following cost function:
    \begin{align*}
        c ((x _1, x _2), (y _1, y _2)) = 
        \begin{cases}
            0, & x _1 = y _1, x _2 = y _2, \\
            101 + |x _1 - y _1|, & x _1 \neq y _1, x _2 = y _2, \\
            100 + |x _1 - y _1| + |x _2 - y _2|, & x _2 \neq y _2.
        \end{cases}
    \end{align*}
    Then $c (x, y) = 0$ iff $x = y$, and $c$ is symmetric.
    Triangular inequality holds because the distance between any two distinct points is between $100$ and $102$. If $\P$ is a uniform distribution on $\{1\} \times [0, 1]$, then $\cM (\hP, \P) = 101$, thus 
    \begin{align*}
        \cLm (101) = \sup _{\P \in \cP (\cX)} \left\lbrace
            \E _{X \sim \P} [f (X)]: \cM  (\hP, \P) \le 101
        \right\rbrace = 1.
    \end{align*}
    However, because $d ((0, t), (1, s)) > 101$ for any $t, s \in [0, 1]$, we would have 
    \begin{align*}
        \E _{\hX \sim \hP} \left[ 
            \sup _{x} \left\lbrace
                f (x): c (\hX, x) \le 101
            \right\rbrace 
        \right] = 0 \neq 1.
        \tag*{\QDEF}
    \end{align*}
\end{remark}

The following result shows that, with additional assumptions on the space and the transport cost function, we can obtain the duality result.
The detailed proofs for Proposition \ref{prop:infty-<} and Theorem \ref{thm:infty-≤} can be found in Appendix \ref{app:maximum-cost}.

\begin{theorem}
    \label{thm:infty-≤}
    Suppose $\cX$ is a Polish space and $c: \cX \times \cX \to [0, \infty)$ is continuous. Let $f$ satisfy Assumption \ref{assum:setup}.
    Then
    \begin{align*}
        \cLm (\rho) &= \E _{\hX \sim \hP} \left[\sup _x \left\lbrace
            f (x): c (\hX, x) \le \rho
        \right\rbrace \right], \\
        (-\cLm) ^* (-\lambda) &= \sup _{\rho \ge 0} \left\lbrace 
        \E _{\hX \sim \hP} \left[\sup _x \left\lbrace
            f (x): c (\hX, x) \le \rho\right] - \lambda \rho
        \right\rbrace
    \right\rbrace .
    \end{align*}
\end{theorem}

\begin{example}[Chance constraint (continued)]
    \label{eg:01-infinity}
    Consider the chance constraint introduced in Example \ref{eg:01} but with $\infty$-Wasserstein distance. 
    If $\cX$ is complete and separable, then the robust chance constraint becomes $\hP (d (\hX, \cSc) \le \rho) \le \beta$.
    In particular, it is infeasible if 
    $\rho \ge d (\supp \hP, \cSc)$.
    Compared with \cite{yang2022wasserstein} which assumes $\cX$ is a normed space, we only require $(\cX, d)$ to be Polish. 
    \QEG
\end{example}

\subsection{Risk-averse Optimization}
\label{sec:risk}

Recall $(\cX, \sF, \hP)$ is a probability space. Consider a concave risk measure $J: \cP (\cX) \to \Rbar$ of the following form:
\begin{align}
    \label{eqn:JP}
    J (\P) := \inf _{\alpha \in A} \E _{X \sim \P} [f _\alpha (X)], 
\end{align}
where $f _\alpha: \cX \to \R$ are a family of measurable functions indexed by $\alpha \in A$, a subset of a linear topological space. 

Given a nominal distribution $\hP$, define 
\begin{align*}
    \cLJ (\rho) & := \sup _{\P \in \cP (\cX)} \left\{
        J (\P): \cK (\hP, \P) \le \rho
    \right\}, \\
    \cLJm (\rho) &:= \sup _{\P \in \cP (\cX)} \left\lbrace
        J (\P): \cM (\hP, \P) \le \rho
    \right\rbrace.
\end{align*}
We assume that there exists a compact set $A' \subset A$ such that for all distributions $\P$ in the distributional uncertainty set, it holds that
\begin{align}
    \label{eqn:achievable}
    \inf _{\alpha \in A} \E _{\P} [f _\alpha] = \min _{\alpha \in A'} \E _{\P} [f _\alpha].
\end{align}
This enables the exchange of sup over $\P$ and inf over $\alpha$ using Sion's minimax theorem.
We have the following result.

\begin{theorem}
    \label{thm:risk-measure}
    
    Let $f _\alpha$ and $c$ satisfy Assumption \ref{assum:setup}, and let $\alpha \mapsto f _\alpha (x)$ be lower semi-continuous and convex for each $x$. 
    Assume for every compact subset $A' \subset A$, $\inf _{\alpha \in A', x \in \cX} f _\alpha (x) > -\infty$. 
    \begin{enumerate}
        \item Suppose $\phi _{\lambda, \alpha} (\hx, x) = f _\alpha (x) - \lambda c (\hx, x)$ satisfies \ref{IP} for every $\lambda > 0$ and $\alpha \in A$. Assume there exists a compact subset $A' \subset A$ such that \eqref{eqn:achievable} holds for every $\P$ satisfying $\cK (\hP, \P) \le \rho$.
        Then 
        \begin{align}
            \label{eqn:cLJ}
            \cLJ (\rho) = \inf _{\lambda \ge 0, \alpha \in A} \left\{
            \lambda \rho + \E _{\hX \sim \hP} \left[ 
                \sup _{x \in \cX} \left\lbrace 
                    f _\alpha (x) - \lambda c (\hX, x)
                \right\rbrace
            \right]\right\}.
        \end{align}

        \item Suppose $\cX$ is a Polish space, and $c: \cX \times \cX \to [0, \infty)$ is continuous. Assume there exists a compact subset $A' \subset A$ such that \eqref{eqn:achievable} holds for every $\P$ satisfying $\cM (\hP, \P) \le \rho$. Then 
        \begin{align}
            \label{eqn:cLJm}
            \cLJm (\rho) = \inf _{\alpha \in A} \E _{\hX \sim \hP} \left[\sup _x \left\lbrace
                f _\alpha (x): c (\hX, x) \le \rho
            \right\rbrace \right].
        \end{align}
    \end{enumerate}
\end{theorem}

Below we list several examples where $X \in \R$ represents a random loss and $A = \R$.

\begin{itemize}
    \item Conditional value-at-risk $\CVaRP _\beta (X)$ at risk level $\beta \in (0, 1)$:
    $
        f _\alpha (x) = \alpha + \frac1{\beta} (x - \alpha) _+
    $.

    \item Variance $\VarP (X)$:
    $
        f _\alpha (x) = (x - \alpha) ^2
    $.

    \item Mean absolute deviation (around median) $\MADP (X)$:
    $
        f _\alpha (x) = |x - \alpha|
    $.

    \item Entropic risk measure $\EntP _\theta (X) = \frac1\theta \log \E [e ^{\theta X}]$ with risk-aversion parameter $\theta > 0$:
    $
        f _\alpha (x) = \alpha + \frac1\theta \left(
            e ^{\theta (x - \alpha)} - 1
        \right)
    $.
\end{itemize}
In Appendix \ref{app:risk}, we can verify \eqref{eqn:achievable} and other assumptions of Theorem \ref{thm:risk-measure} hold for all these risk measures.

\begin{example}
\label{eg:risk}
    Let $Z \in (\R ^d, \|\cdot\|)$ be the vector of random loss of $d$ assets with nominal distribution $\hQ$, and $b \in \R ^d$ be a portfolio weight vector on these assets.
    The portfolio loss is thus $b ^\top Z$.
    We have the following results on the robust risk of portfolio loss under various risk measures. Their proofs are given in Appendix \ref{app:risk}.

    \begin{itemize}
        \item \label{eg:cvar} 
        Conditional value-at-risk: for $p \in [1, \infty]$:
        \begin{align*}
            \sup _{\Q \in \cP (\cZ)} \left\lbrace
                \CVaRQ _\beta (b ^\top Z): \cW _p (\hQ, \Q) \le \rho
            \right\rbrace = \CVaRhQ _\beta (b ^\top \hZ) + (1 - \beta) ^{-\frac1p} \| b \| _* \rho.
        \end{align*} 
        
        \item \label{eg:var}
        Variance: for $p = 2$,    
        \begin{align*}
            \sup _{\Q \in \cP (\cZ)} \left\lbrace
                \VarQ (b ^\top Z): \cW _2 (\hQ, \Q) \le \rho
            \right\rbrace = \left(
                \VarhQ (b ^\top \hZ) ^\frac12 + \| b \| _* \rho
            \right) ^2.
        \end{align*}
        For $p = \infty$, 
        \begin{align*}
            \sup _{\Q \in \cP (\cZ)} \left\lbrace
                \VarQ (b ^\top Z): \cW _\infty (\hQ, \Q) \le \rho
            \right\rbrace = \min _{\alpha \in \R} \E _{\hZ \sim \hQ} \left[
                (|\hZ - \alpha| + \| b \| _* \rho) ^2
            \right].
        \end{align*}
        For $1 \le p < 2$, the robust loss is positive infinity.
        
        \item \label{eg:mad}
        Mean average deviation:
        for $1 \le p \le \infty$, 
        \begin{align*}
            \sup _{\Q \in \cP (\cZ)} \left\lbrace
                \MADQ (b ^\top Z): \cW _p (\hQ, \Q) \le \rho
            \right\rbrace = \MADhQ (b ^\top \hZ) + \| b \| _* \rho.
        \end{align*}
        
        \item \label{eg:entropic}
        Entropic risk measure:
        for $p = \infty$,
        \begin{align*}
            \sup _{\Q \in \cP (\cZ)} \left\lbrace
                 \EntQ _\theta (b ^\top Z): \cW _\infty (\hQ, \Q) \le \rho
            \right\rbrace = \EnthQ _\theta (b ^\top \hZ)] + \| b \| _* \rho.
        \end{align*}
        For $1 \le p < \infty$, the robust loss is positive infinity. 
        \QEG
    \end{itemize}
\end{example}

\subsection{Globalized Distributionally Robust Counterpart}
\label{sec:globalized}

The globalized distributionally robust counterpart \cite{liu2023globalized} studies the following problem
\begin{align}
    \label{eqn:G}
    \tag{\textsf{G}}
    \sup _{\P, \tP \in \cP (\cX)} \left\{ 
        \E _{X \sim \P} [f (X)] - \lambda \cK (\tP, \P): \tcK (\hP, \tP) \le \theta
    \right\}.
\end{align}
We also consider its hard- and soft-constrained variants
\begin{align}
    \label{eqn:G-hard}
    \tag{\textsf{G-hard}}
    \cLG (\rho, \theta) := 
    \sup _{\P, \tP \in \cP (\cX)} \left\{ 
        \E _{X \sim \P} [f (X)]: \cK (\tP, \P) \le \rho, \tcK (\hP, \tP) \le \theta
    \right\}, \\
    \label{eqn:G-soft}
    \tag{\textsf{G-soft}}
    \sup _{\P, \tP \in \cP (\cX)} \left\{ 
        \E _{X \sim \P} \left[
            f (X)
        \right] - \lambda \cK (\tP, \P) - \mu \tcK (\hP, \tP)
    \right\}.
\end{align}

The following result extends the work of \cite{liu2023globalized}, which is based on the assumption that $\cX$ is a subset of Euclidean space and the transport cost is defined by the 1-Wasserstein distance. The proof can be found in Appendix \ref{app:globalized}.

\begin{proposition}
    \label{prop:globalized}
    Let the loss function $f$ and two cost functions $c, \tc$ satisfy Assumption \ref{assum:setup}.
    For $\rho, \theta > 0$, $\lambda, \mu \ge 0$, if $(\tx, x) \mapsto f (x) - \lambda c (\tx, x)$ satisfies \ref{IP} for every $\lambda \ge 0$, and $(\hx, x) \mapsto \sup _{\tx \in \cX} f (x) - \lambda c (\tx, x) - \mu \tc (\hx, \tx)$ satisfies \ref{IP} for every $\lambda, \mu \ge 0$, then for every $\lambda, \mu, \rho, \theta > 0$,
    \eqref{eqn:G-hard} is equivalent to 
    \begin{align*}
        \min _{\lambda, \mu \ge 0} \left\{
            \lambda \rho + \mu \theta + \E _{\hX \sim \hP} \left[
                \sup _{x, \tx \in \cX} \left\{
                    f (x) - \lambda c (\tx, x) - \mu \tc (\hX, \tx)
                \right\}
            \right]
        \right\}, 
    \end{align*}
    \eqref{eqn:G} is equivalent to 
    \begin{align*}
        (-\cLG (\cdot, \theta)) ^* (-\lambda) &= \min _{\mu \ge 0} \left\{
            \mu \theta + \E _{\hX \sim \hP} \left[
                \sup _{x, \tx \in \cX} \left\{
                    f (x) - \lambda c (\tx, x) - \mu \tc (\hX, \tx)
                \right\}
            \right]
        \right\}, 
    \end{align*}
    and \eqref{eqn:G-soft} is equivalent to 
    \begin{align*}
        \cLG \doublestar (-\lambda, -\mu) &= 
        \E _{\hX \sim \hP} \left[
            \sup _{x, \tx \in \cX} \left\{
                f (x) - \lambda c (\tx, x) - \mu \tc (\hX, \tx)
            \right\}
        \right].
    \end{align*}
    Here $\cLG \doublestar (-\lambda, \cdot)$ is the dual of the mapping $\theta \mapsto -(-\cLG (\cdot, \theta)) ^* (-\lambda)$.
    Moreover, if $c  = \tc  = d $, then \eqref{eqn:G} is further equivalent to 
    \begin{align*}
        \min _{0 \le \mu \le \lambda} \left\{
            \mu \theta + \E _{\hX \sim \hP} \left[
                \sup _{x \in \cX} \left\{
                    f (x) - \mu d (\hX, x)
                \right\}
            \right]
        \right\}.
    \end{align*}
\end{proposition}

\section{Concluding Remarks}\label{sec:conclusion}
We have developed a new duality proof for Wasserstein distributionally robust optimization. 
The new result offers an alternative view of duality from the perspective of the interchangeability principle. This suggests that establishing duality in broader settings hinges primarily on verifying the interchangeability property, which has been more extensively explored.

\bigskip

\begin{APPENDICES}
\SingleSpacedXI

\section{Auxiliary Results}\label{app:aux}

\begin{lemma}
    \label{lem:convex}
    Assume Assumption \ref{assum:setup} holds. Recall for $\lambda \ge 0$, 
    \begin{align*}
        (-\cL) ^* (-\lambda) &= \sup _{\P \in \cPs} \left\lbrace
            \E _{X \sim \P} [ f (X) ] - \lambda \cK (\hP, \P)
        \right\rbrace, 
        &
        \cG (\lambda) &= \E _{\hX \sim \hP} \left[ 
            \sup _{x \in \cX} \left\lbrace 
                f (x) - \lambda c (\hX, x) 
            \right\rbrace 
        \right].
    \end{align*}
    Then $(-\cL) ^* (-\cdot)$ and $\cG (\cdot)$ are lower bounded by $\E_{\hP} [f]$, monotonically decreasing, convex, and lower semi-continuous on $[0, \infty)$. 
\end{lemma}

\begin{lemma}
    \label{lem:convex-2}
    If $f: \R \to \R \cup \{+\infty\}$ is a monotonically decreasing convex function with $f (\rho) = +\infty$ for $\rho < 0$ and $f \not\equiv +\infty$, then $f ^* (-\lambda) = \sup _{\rho \in \R} \{-\lambda \rho - f (\rho)\}$ is a lower semi-continuous, monotonically decreasing convex function of $\lambda$ with $f ^* (-\lambda) = +\infty$ for $\lambda < 0$ and $f ^* \not\equiv +\infty$.    
\end{lemma}

\begin{lemma}
    \label{lem:f-lc}
    $\phi: \cX \times \cX \to \Rninfty$ is diagonally dominant if and only if there exists a measurable function $f: \cX \to \Rninfty$, a constant $\lambda \ge 0$, and a measurable function $c: \cX \times \cX \to [0, \infty]$ vanishing on diagonal such that $\phi (\hx, x) = f (x) - \lambda c (\hx, x)$.
\end{lemma}

\end{APPENDICES}

\bibliographystyle{informs2014}
\bibliography{ref.bib}

\ECSwitch

\ECHead{\centering Additional Proofs} \label{app:eg}

\section{Proofs for Appendix \ref{app:aux}}

\proof{Proof of Lemma \ref{lem:concave}.}
    The monotonicity of $\RP (\rho)$ can be seen from the definition. Moreover, since $\cK (\hP, \hP) = 0$, 
    \begin{align*}
        \RP (\rho) \ge \RP (0) \ge \E _{X \sim \hP} [f (X)] > -\infty.
    \end{align*}
    Therefore for all $\rho \ge 0$, $\RP (\rho)$ is bounded from below.
    To verify the concavity, fix $\rho _0, \rho _1 \ge 0$. Pick any $t\in[0,1]$ and $\P _0, \P _1 \in \cP (\cX)$ satisfying $\cK (\hP, \P _j) \le \rho _j$ and $\E _{\P_j} [f] > -\infty$, $j = 0, 1$, and denote $\P _t = (1 - t) \P _0 + t \P _1$. For arbitrary $\epsilon > 0$, we can find transport plans $\gamma _0 \in \Gamma (\hP, \P _0)$, $\gamma _1 \in \Gamma (\hP, \P _1)$ such that $\E _{\gamma _0} [c] \le \cK (\hP, \P _0) + \epsilon$, $\E _{\gamma _1} [c] \le \cK (\hP, \P _1)  + \epsilon$. Define $\gamma _t = (1 - t) \gamma _0 + t \gamma _1$, then $\gamma _t \in \Gamma (\hP, \P _t)$ and 
    \begin{align*}
        \cK (\hP, \P _t) \le \E _{\gamma _t} [c] = (1 - t) \E _{\gamma _0} [c] + t \E _{\gamma _1} [c] \le (1 - t) (\cK (\hP, \P _0) + \epsilon) + t (\cK (\hP, \P _1) + \epsilon) \\
        \le (1 - t) \cK (\hP, \P _0) + t \cK (\hP, \P _1) + \epsilon \le (1 - t) \rho _0 + t \rho _1 + \epsilon.
    \end{align*}
    Since it is true for any $\epsilon$, we know $\cK (\hP, \P _t) \le (1 - t) \rho _0 + t \rho _1$,
    hence $\P_t$ is a feasible solution to \eqref{problem}  {with $\rho = (1 - t) \rho _0 + t \rho _1$} and
    \begin{align*}
        \RP ((1 - t) \rho _0 + t \rho _1) 
        &\ge 
        \E _{X \sim \P _t} [f (X)] %
        = 
        (1 - t)  
            \E _{X \sim \P _0} [f (X)] %
         + t 
            \E _{X \sim \P _1} [f (X)] %
        .
    \end{align*}
    Taking the supremum over $\P _0$ and $\P _1$, we have
    \begin{align*}
        \RP ((1 - t) \rho _0 + t \rho _1) \ge (1 - t) \RP (\rho _0) + t \RP (\rho _1),
    \end{align*}
    which completes the proof. 
    \QED
\endproof

\proof{Proof of Lemma \ref{lem:convex-2}.}
    Since $f (\rho) = +\infty$ for $\rho < 0$, we have 
    \begin{align*}
        f ^* (-\lambda) = \sup _{\rho \in \R} \{-\lambda \rho - f (\rho)\} = \sup _{\rho \ge 0} \{-\lambda \rho - f (\rho)\}.
    \end{align*}
    For each fixed $\rho \ge 0$, $-\lambda \rho - f (\rho)$ is a monotonically decreasing, lower semi-continuous convex function of $\lambda$, so the supremum over $\rho$ is also monotonically decreasing, lower semi-continuous, and convex. Suppose $f (\rho _0) < +\infty$ at some $\rho _0 \ge 0$. For each $\lambda < 0$, 
    \begin{align*}
        f ^* (-\lambda) = \sup _{\rho \ge 0} \{-\lambda \rho - f (\rho)\} \ge \sup _{\rho \ge \rho _0} \{-\lambda \rho - f (\rho)\} \ge \sup _{\rho \ge \rho _1} \{-\lambda \rho - f (\rho _0)\} = +\infty.
    \end{align*}
    Pick $\rho _1 > \rho _0$.
    For each $\lambda \ge 0$,
    \begin{align*}
        f ^* (-\lambda) = \sup _{\rho \ge 0} \{-\lambda \rho - f (\rho)\} = \sup _{0 \le \rho \le \rho _1} \{-\lambda \rho - f (\rho)\} \vee \sup _{\rho \ge \rho _1} \{-\lambda \rho - f (\rho)\}.
    \end{align*}
    When $0 \le \rho \le \rho _1$, $-\lambda \rho - f (\rho) \le - f (\rho _1) < +\infty$. When $\rho \ge \rho _1$, by convexity we have $f (\rho) \ge f (\rho _1) + (\rho - \rho _1) \frac{f (\rho _1) - f (\rho _0)}{\rho _1 - \rho _0}$, so if $\lambda \ge \frac{f (\rho _0) - f (\rho _1)}{\rho _1 - \rho _0}$ we must have 
    \begin{align*}
        -\lambda \rho - f (\rho) \le -\lambda \rho - f (\rho _1) - (\rho - \rho _1) \frac{f (\rho _1) - f (\rho _0)}{\rho _1 - \rho _0} \le - f (\rho _1) - \rho _1 \frac{f (\rho _0) - f (\rho _1)}{\rho _1 - \rho _0} < -f (\rho _1) < +\infty.
    \end{align*}
    Hence $f ^* (\lambda) \le -f (\rho _1) < +\infty$, so $f ^* \not\equiv +\infty$.
    \QED
\endproof

\proof{Proof of Lemma \ref{lem:convex}.}
    The lower bound follows by setting $x = \hX$: for any $\lambda \in [0, \infty)$, we have
    \begin{align*}
        (-\cL) ^* (-\lambda) &\ge \E _{\hX \sim \hP} [f (\hX)] - \lambda \cK (\hP, \hP) = \E _{\hP} [f], 
        &
        \cG ^* (\lambda) &\ge \E _{\hX \sim \hP} \left[f (\hX) - \lambda c (\hX, \hX)\right] = \E _{\hP} [f].
    \end{align*}
    Here we used $\cK (\hP, \hP) = 0$ because $c (x, x) = 0$ for every $x \in \cX$.

    If $\cL (\rho) < +\infty$ for all $\rho > 0$, then
    $(-\cL) ^* (-\cdot)$ is decreasing, convex, and lower semi-continuous because of Lemma \ref{lem:convex-2}. If $\cL (\rho) = +\infty$ for all $\rho > 0$, then $(-\cL) ^* \equiv +\infty$.
    
    Since $f (x) - \lambda c (\hx, x)$ is decreasing and affine in $\lambda$, $\Phi (\lambda; \hx) := \sup _{x \in \cX} \left\lbrace 
        f (x) - \lambda c (\hx, x) 
    \right\rbrace$ is also a decreasing, convex, and lower semi-continuous function of $\lambda$. Recall that when $\lambda = 0$ and $c (\hx, x) = +\infty$, we use the convention $0 \cdot \infty = \infty$. We now verify that 
    \begin{enumerate}
        \item[(a)] $\cG ^*$ is monotonically decreasing: 
        \begin{align*}
            \lambda _1 \le \lambda _2 \implies \Phi (\lambda _1; \hx) \ge \Phi (\lambda _2; \hx) \text{ for all } \hx \in \cX \implies \E _{\hP} [\Phi (\lambda _1; \hX)] \ge \E _{\hP} [\Phi (\lambda _2; \hX)].
        \end{align*}

        \item[(b)] $\cG ^*$ is convex: 
        \begin{align*}
            \lambda _\theta = (1 - \theta) \lambda _1 + \theta \lambda _2 
            &\implies \Phi (\lambda _\theta; \hx) \le (1 - \theta) \Phi (\lambda _0; \hx) + \theta \Phi (\lambda _1; \hx) \text{ for all } \hx \in \cX \\
            &\implies \E _{\hP} [\Phi (\lambda _\theta; \hX)] \le (1 - \theta) \E _{\hP} [\Phi (\lambda _0; \hX)] + \theta \E _{\hP} [\Phi (\lambda _1; \hX)].
        \end{align*}

        \item[(c)] $\cG ^*$ is lower semi-continuous: note that $\Phi (\lambda; \hx) \ge f (\hx)$ and $\E _{\hP} [f] > -\infty$. Taking $\lambda _n \to \lambda$ where $\lambda _n, \lambda \in [0, \infty)$, then by Fatou's lemma, we have
        \begin{align*}
            \liminf _{\lambda _n \to \lambda} \E _{\hP} [\Phi (\lambda _n; \hX)] \ge \E _{\hP} \left[\liminf _{\lambda _n \to \lambda} \Phi (\lambda _n; \hX)\right] \ge \E _{\hP} [\Phi (\lambda)].
        \end{align*}
    \end{enumerate}
    With this we complete the proof.
\QED
\endproof

\begin{proof}{Proof of Lemma \ref{lem:f-lc}.}
    $f (x) - \lambda c (\hx, x) \le f (x) - 0 = f (x) - \lambda c (x, x)$, so $f - \lambda c$ is diagonally dominant. If $\phi$ is diagonally dominant, we define $\lambda := 1$, $f (x) := \phi (x, x)$, and $c (\hx, x) := f (x) - \phi (\hx, x)$ when $f (x) > -\infty$, $c (\hx, x) = 0$ when $f (x) = -\infty$. Then $c (\hx, x) \ge 0$ and $c (x, x) = 0$.
\end{proof}

\section{Proof of Proposition \ref{prop:opt-selection}}

Before proving Proposition \ref{prop:opt-selection}, we make a simple observation. 

\begin{lemma}
    $A \subset \cX \times \cX$ is a diagonally dominant set if and only if its indicator function $\mathbf1 _A$ is a diagonally dominant function. $\phi: \cX \times \cX \to \Rninfty$ is a diagonally dominant function if and only if its superlevel set $\{\phi > \alpha\}$ is a diagonally dominant set for any $\alpha \in \R$. $E: \cX \to \sF \setminus \lbrace\varnothing\rbrace$ is a diagonally dominant set-valued function if and only if its graph is a diagonally dominant set.
\end{lemma}

\proof{Proof of Proposition \ref{prop:opt-selection}.}
    We first prove the sufficiency. Let $\phi$ be a diagonally dominant $(\sF\otimes\sF)$-measurable function.
    Define $\Phi (\hx) = \sup _{x \in \cX} \phi (\hx, x)$.
    For any $\alpha \in \R$, the superlevel set of $\Phi$ can be regarded as
    \begin{align*}
        \left\lbrace \hx: \Phi (\hx) > \alpha \right\rbrace = 
        \left\lbrace \hx: \exists x, \phi (\hx, x) > \alpha \right\rbrace =
        \Proj _{\hx} (\left\lbrace (\hx, x): \phi (\hx, x) > \alpha \right\rbrace).
    \end{align*}
    The superlevel set $\{\phi > \alpha\}$ is diagonally dominant.
    By assumption \ref{ass:Proj}, $\Proj _{\hx}$ maps $(\sF \otimes \sF)$-measurable diagonally dominant sets to $\sF_{\hP}$-measurable sets, thus the superlevel set of $\Phi$ is $\sF_{\hP}$-measurable. Therefore, $\Phi$ is $\hP$-measurable, and it remains to show that 
    \begin{align*}
        \E _{\hP} [\Phi] = \sup _{\gamma \in \Gamma _{\hP}} \E _{\gamma} [\phi].
    \end{align*}
    Since for any $x \in \cX$, $\phi (\hx, x) \le \Phi (\hx)$, it is clear that for any $\gamma \in \Gamma _{\hP}$,
    \begin{align*}
        \E _{\hX \sim \hP} \left[
            \Phi (\hX)
        \right] \ge \E _{(\hX, X) \sim \gamma} \left[
            \phi (\hX, X)
        \right].
    \end{align*}
    To see the other direction, we may assume $\E  _{\hX \sim \hP} [\Phi (\hX)] > -\infty$, otherwise the conclusion holds trivially. Then $\{\Phi = -\infty\}$ is a $\hP$-nullset. We fix $\epsilon, M > 0$, and below we construct a near optimal $\gamma \in \Gamma _{\hP}$.
    
    Define $S _n = \{\hx \in \cX: n \epsilon < \Phi (\hx) \le (n + 1) \epsilon \}$ for $n \in \mathbb Z$. $S _n$ is $\sF _{\hP}$-measurable, so we can find $B _n \subset S _n$ which is $\sF$-measurable and $\hP (S _n \setminus B _n) = 0$.
    Define set-valued function $E _n: \cX \to \sF \setminus \{\varnothing\}$ by
    \begin{align*}
        E _n (\hx) = \begin{cases}
            \cX \setminus B _n & \hx \notin B _n \\
            \{x \in \cX:  \phi (\hx, x) > n \epsilon \} & \hx \in B _n
        \end{cases}.
    \end{align*}
    $\Graph (E _n) = (\cX \setminus B _n) \times (\cX \setminus B _n) \cup ((B _n \times \cX) \cap \{\phi > n \epsilon\})$ is $(\sF \otimes \sF)$-measurable.
    We claim $E _n$ is diagonally dominant. That is, $x \in E _n (\hx)$ implies $x \in E _n (x)$. To see this, note that if $x \notin B _n$ then $x \in \cX \setminus B _n = E _n (x)$; if $\hx \notin B _n$ and $x \in E _n (\hx) = \cX \setminus B _n$ then $x \notin B _n$ so $x \in E _n (x)$. Now suppose $\hx, x \in B _n$ and $x \in E _n (\hx)$, then $\phi (x, x) \ge \phi (\hx, x) > n \epsilon$, so again we have $x \in E _n (x)$. This finishes the proof of the claim. By assumption \ref{ass:Sel*} we can find a measurable transport plan $\gamma _n \in \Gamma _{\hP}$ supported in $\Graph (E _n)$.

    Define $S _\infty = \{\hx \in \cX: \Phi (\hx) = \infty\}$. $S _\infty$ is $\sF _{\hP}$-measurable, so we can find $B _\infty \subset S _\infty$ which is $\sF$-measurable and $\hP (S _\infty \setminus B _\infty) = 0$.
    For some $M > 0$ to be determined, define set-valued function $E _\infty: \cX \to \sF \setminus \{\varnothing\}$ by
    \begin{align*}
        E _\infty (\hx) = \begin{cases}
            \cX \setminus B _\infty & \hx \notin B _\infty \\
            \{x \in \cX:  \phi (\hx, x) > M \} & \hx \in B _\infty
        \end{cases}.
    \end{align*}
    Same as before, $E _\infty$ is diagonally dominant, so we can find a measurable transport plan $\gamma _\infty \in \Gamma _{\hP}$ supported in $\Graph (E _\infty)$.

    Now we define measure $\gamma \in \cP (\cX \times \cX, \sF \otimes \sF)$ by
    \begin{align*}
        \gamma (A) = \sum _{n \in \mathbb Z \cup \{+\infty\}} \gamma _n (A \cap (B _n \times \cX)).
    \end{align*}
    Then for any $S \subset \cX$, 
    \begin{align*}
        \gamma (S \times \cX) = \sum _{n \in \mathbb Z \cup \{+\infty\}} \gamma _n ((S \cap B _n) \times \cX) = \sum _{n \in \mathbb Z \cup \{+\infty\}} \hP (S \cap B _n) = \hP (S),
    \end{align*}
    so $\gamma \in \Gamma _{\hP}$.
    In the last equality, we used that $B _n \subset S _n$ are pairwise disjoint and $\hP (\cX \setminus \bigcup _{n \in \mathbb Z \cup \{+\infty\}} B _n) = 0$. Moreover, $\gamma$ is supported in 
    \begin{align*}
        \{(\hx, x) \in \cX \times \cX: \phi (\hx, x) > \Phi (\hx) - \epsilon \text{ if } \Phi (\hx) < \infty, \phi (\hx, x) > M \text{ if } \Phi (\hx) = \infty \}.
    \end{align*}
    Therefore, $\E _{\gamma} [\phi] \ge \E _{\hP} [(\Phi - \epsilon) \mathbf1 \{\Phi < +\infty\}] + M \hP (\Phi = +\infty)$. By making $\epsilon$ arbitrarily small and $M$ arbitrarily large, we have $\sup _{\gamma \in \Gamma _{\hP}} \E _{\gamma} [\phi] \ge \E _{\hP} [\Phi]$. This proves that \ref{ass:Proj} and \ref{ass:Sel*} combined imply \ref{IP} holds for all diagonally dominant functions.
    
    Next, we prove the necessity. Suppose \ref{IP} holds for all diagonally dominant functions.  
    Given $A \in \sF\otimes\sF$ diagonally dominant, let $\phi$ be the indicator of the set $A$: $\phi(\hx,x)=1$ if $(\hx,x)\in A$ and $0$ otherwise. Then $\phi$ is $\sF\otimes\sF$-measurable, diagonally dominant, and by \ref{IP} the function $\Phi(\hx)= \sup _{x\in \cX} \phi(\hx, x)$ is $\hP$-measurable. Observe that
    \begin{align*}
        \Proj _{\hx} (A) = \{\hx\in\cX: \Phi(\hx) \ge 1\},
    \end{align*}
    which is the upper level set of $\Phi$, and thus belongs to $\sF_{\hP}$. Therefore \ref{IP} implies \ref{ass:Proj}. 
    Lastly, given a diagonally dominant set function $E: \cX \to \sF \setminus \lbrace\varnothing\rbrace$, let
    \begin{align*}
        \phi (\hx, x) = \begin{cases}
            0 & x \in E (\hx) \\
            -\infty & x \notin E (\hx)
        \end{cases}.
    \end{align*}
    That is, $\phi = -\infty \cdot  \mathbf 1 _{\cX \setminus \Graph (E)}$. To see it is diagonally dominant, note that 
    \begin{align*}
        \phi (\hx, x) = 0 \implies x \in E (\hx) \implies x \in E (x) \implies \phi (x, x) = 0.
    \end{align*}
    So $\phi (\hx, x) \le \phi (x, x)$. Since $E (\hx) \neq \varnothing$, $\Phi (\hx) := \sup _{x \in \cX} \phi (\hx, x) = 0$. By \ref{IP}, 
    \[
        0 = \E _{\hX \sim \hP} \left[ 
            \Phi (\hX, x) 
        \right]
        = 
        \sup _{\gamma \in \Gamma _{\hP}} \left\lbrace
            \E _{(\hX, X) \sim \gamma} [ \phi (\hX, X) ]
        \right\rbrace.
    \]
    Note that 
    \[
        \E _{\gamma} [\phi] = \begin{cases}
            0, & \supp \gamma \subset \Graph (E) ,\\
            -\infty, & \text{otherwise}.
        \end{cases}
    \]
    Hence, there exists some $\gamma \in \Gamma _{\hP}$ supported in $A$. Therefore \ref{IP} implies \ref{ass:Sel*}.
\QED
\endproof

\proof{Proof of Proposition \ref{prop:c-continuous}.}
    Note that $\phi$ is continuous in $\hx$, so
    $\Phi (\hx) = \sup _{x \in \cX} \phi (\hx, x)$ is lower semi-continuous, and thus Borel measurable. Therefore, for any $\gamma \in \Gamma _{\hP}$, it holds that $\E _\gamma [\phi] \le \E _{\hP} [\Phi]$. It remains to construct an $\epsilon$-optimizer $\gamma$. 
    
    \newcommand{\hK}{\widehat K}
    
    First, we assume $\E _{\hP} [\Phi] < +\infty$, and fix $\epsilon > 0$. Since $\hP$ is tight, there exists a compact subset $\hK \subset \subset \cX$ sufficiently large such that 
    $\E _{\hP} [|f| \mathbf1 _{\hK ^c}] < \epsilon$ and
    $\E _{\hP} [|\Phi| \mathbf1 _{\hK ^c}] < \epsilon$. 
    Moreover, fix some $\hx _0 \in \cX$, we define 
    \begin{align*}
        K _n = \hK \cup B _n (\hx _0) = \hK \cup \{x \in \cX: d (\hx _0, x) \le n\},
    \end{align*}
    and define 
    \begin{align*}
        \phi _n (\hx, x) = \phi (\hx, x) - \infty \mathbf1\{x \notin K _n\}, \qquad 
        \Phi _n (\hx) = \sup _{x \in K _n} \phi (\hx, x) = \sup _{x \in \cX} \phi _n (\hx, x).
    \end{align*}
    Then $\Phi _n \to \Phi$ pointwise as $n \to \infty$. Since $f \le \Phi _n \le \Phi$ in $\hK$, From dominant convergence theorem, we have 
    \begin{align*}
        \E _{\hP} [\Phi \mathbf1 _{\hK}] - \E _{\hP} [\Phi _n \mathbf1 _{\hK}] < \epsilon    
    \end{align*}
    for $n$ sufficiently large. We fix $n$ from now on.
    Note that $\phi (\hx, x)$ is uniformly continuous in $\hx$ in $\hK \times K _n$: there exists $\delta > 0$ such that if $\hx _1, \hx _2 \in \hK$ and $d(\hx _1, \hx _2) < \delta$, we must have $|\phi (\hx _1, x) - \phi (\hx _2, x)| \le \epsilon$ for all $x \in K _n$, and consequently $|\Phi _n (\hx _1) - \Phi _n (\hx _2)| \le \epsilon$.
    Since $\hK$ is compact, there exists a $\delta$-net $\widehat{\cX} = \{ \hx _i \} _{i = 1} ^n \subset \hK$. 
    Define $U _i = \hK \cap B _\delta (\hx _i) \setminus \bigcup _{j < i} B _\delta (\hx _j)$, then $\{U _i\} _{i = 1} ^n$ forms a partition of $\hK$. 
    For each $\hx _i$, we can find $x _i$ such that 
    \begin{align*}
        \phi (\hx _i, x _i) > \Phi (\hx _i) - \epsilon.
    \end{align*}
    Now we construct a Borel-measurable selection mapping 
    \begin{align*}
        T (\hx) = \begin{cases}
            x _i & \hx \in U _i \\
            x & \hx \in \hK ^c
        \end{cases}, \qquad \hx \in \cX.
    \end{align*}
    This induces a measure $\gamma = (\Id \otimes T) _\# \hP$. Under this selection, we have 
    \begin{align*}
        \E _{\gamma} [\phi] = \E _{\hP} \left[
            \phi (\hX, T (\hX))
        \right] = \E _{\hP} \left[
            \phi (\hX, \hX) \mathbf1 \{\hX \in \hK ^c\}
        \right] + \sum _{i = 1} ^n \E _{\hP} \left[
            \phi (\hX, x _i) \mathbf1 \{\hX \in U _i\}
        \right].
    \end{align*}
    The first term is 
    \begin{align*}
        \E _{\hP} \left[
            \phi (\hX, \hX) \mathbf1 \{\hX \in \hK ^c\}
        \right] = \E _{\hP} [f \mathbf1 _{\hK ^c}] \ge -\epsilon.
    \end{align*}
    For the second term, note that $|\hx - \hx _i| < \delta$ for $\hx \in U _i$, so by uniform continuity we have 
    \begin{align*}
        \phi (\hx, x _i) \ge \phi _n (\hx, x _i) \ge \phi _n (\hx _i, x _i) - \epsilon \ge \Phi _n (\hx _i) - 2 \epsilon \ge \Phi _n (\hx) - 3 \epsilon.
    \end{align*}
    Hence we have 
    \begin{align*}
        \E _{\gamma} [\phi] \ge -4\epsilon + \sum _{i = 1} ^n \E _{\hP} \left[
            \Phi _n (\hX) \mathbf1 \{\hX \in U _i\}
        \right] = -4\epsilon + \E _{\hP} [\Phi _n \mathbf1 _{\hK}] > -5\epsilon + \E _{\hP} [\Phi \mathbf1 _{\hK}] > -6\epsilon + \E _{\hP} [\Phi].
    \end{align*}
    Since $\epsilon$ can be chosen arbitrarily small, we proved interchangeability for $\phi$. The case $\E _{\hP} [\Phi] = +\infty$ is similar and we omit the proof.
\QED
\endproof

\begin{remark}
    \label{rmk:c-continuous}
    It can be seen from the proof that beyond the Wasserstein setting, we need $c$ to be continuous in $\hx$ uniformly in $\hK \times K _n$ for any compact $\hK$ and some sequence of $K _n \supset \hK$ such that $\bigcup _n K _n = \cX$. In particular, this would hold if $c$ is continuous and $\cX$ is a $\sigma$-compact metrizable topological space.
\end{remark}

\section{Proofs for Section \ref{eg:generality}}
\label{app:generality}

In this section, we provide additional details of the proof in Example \ref{eg:MDP}, Example \ref{eg:SP}, and Example \ref{eg:01}.

\proof{Proof of Example \ref{eg:MDP}.}
    We start with $t=T$. We have
    \begin{align*}
        V _T (s) = \inf _{a \in \cA (s)} g _T (s, a),    
    \end{align*}
    which is lower semi-analytic \citep[Proposition 7.47]{bertsekas1996stochastic}, and in particular, $\hP$-measurable \citep[Corollary 7.42.1]{bertsekas1996stochastic}.
    Since $g_T$ is bounded from below by our assumption, $V_T$ is also bounded from below. In a Borel space or Polish space, any measure is tight.
    Thus by Proposition \ref{prop:c-continuous} we have
    \begin{align*}
        V _{T - 1} (s) & = \inf _{a \in \cA (s)} \left\{ 
            g _{T - 1} (s, a) + \sup _{\P \in \frakM (s, a)} \E _{s' \sim \P (\cdot \mid s, a)} [V _{T} (s')] 
        \right\} \\
        & = \inf_{\substack{a \in \cA(s) \\ \lambda \ge 0}} \left\{ 
            g _{T - 1} (s, a) + \lambda \rho (s, a) + \E _{\hs' \sim \hP (\cdot \mid s, a)} \left[
                \sup _{s' \in \cS} \big\{ 
                    V _{T} (s') - \lambda c (\hs', s')
                \big\} 
            \right] 
        \right\},
    \end{align*}
    which is also bounded from below since $g _{T - 1}$ is bounded from below by assumption and $\rho > 0$.
  
    Suppose we have shown $V _{t + 1} (\cdot)$ is lower semi-analytic, $V _{t}$ is bounded from below and obtain the reformulation for $V _{t}$. Now let us show that $V _t (\cdot)$ is lower semi-analytic and derive the expression for $V _{t-1}$ and show it is bounded from below.
    By the continuity of $c$, $\hs' \mapsto c (\hs', s')$ is continuous for each $s'$, so the function $\hs' \mapsto \sup _{s' \in \cS} \big\{ 
        V _{t + 1} (s') - \lambda c (\hs', s')
    \big\} $ is the supremum of a family of continuous function, which is lower semi-continuous, and furthermore is Borel measurable and thus lower semi-analytic. 
    Hence the function $(s, a) \mapsto \E _{\hs' \sim \hP (\cdot \mid s, a)} \left[
        \sup _{s' \in \cS} \big\{ 
            V _{t + 1} (s') - \lambda c (\hs', s')
        \big\} 
    \right] $ is lower semi-analytic \citep[Proposition 7.48]{bertsekas1996stochastic}.
    Since $g_{t-1}$ and $\rho$ are lower semi-analytic due to our assumptions, $V_t$ is also lower semi-analytic \citep[Proposition 7.47]{bertsekas1996stochastic}. Then using Proposition \ref{prop:c-continuous} again we have
    \begin{align*}
        V _{t - 1} (s) & = \inf _{a \in \cA (s)} \left\{ 
            g _{t - 1} (s, a) + \sup _{\P \in \frakM (s, a)} \E _{s' \sim \P (\cdot \mid s, a)} [V _{t} (s')] 
        \right\} \\
        & = \inf_{\substack{a \in \cA(s) \\ \lambda \ge 0}} \left\{ 
            g _{t - 1} (s, a) + \lambda \rho (s, a) + \E _{\hs' \sim \hP (\cdot \mid s, a)} \left[
                \sup _{s' \in \cS} \big\{ 
                    V _{t} (s') - \lambda c (\hs', s')
                \big\} 
            \right] 
        \right\},
    \end{align*}
    which is bounded from below since $g _{t - 1}, \rho, V _t$ are all bounded from below. Therefore the proof is completed.
    \QED
  \endproof

\proof{Proof of Example \ref{eg:SP}.}
    We start with $t=T$. In this case,
    \[
      Q_T(u_{T-1},x_T) = \inf_{u\in\cU_T(u_{T-1},x_T)} f_T(u,x_T).
    \]
    Since $f_T$ is random lower semi-continuous and $\cU_T$ is uniformly bounded, $Q_{T}(\cdot,\cdot)$ is random lower semi-continuous \citep[Theorem 9.50]{shapiro2021lectures}, and particularly, $Q _{T} (u _{T - 1}, \cdot)$ is measurable.
    Since $f_T$ is bounded from below by our assumption, $Q_T$ is also bounded from below.
    Thus Assumption \ref{assum:setup} holds and using Example \ref{eg:suslin} we have
    \begin{align*}
        Q _{T - 1} (u _{T - 2}, x _{T - 1}) &= \inf _{u \in \cU _{T - 1} (u _{T - 2}, x _{T - 1})} \left\{ 
            f _{T - 1} (u, x _{T - 1}) + \sup _{\P \in \frakM _{T}}
                \E _{\P} \big[
                    Q _{T} (u, x _{T})
                \big]
            \right\} \\
        & = \inf _{\substack{u \in \cU _{T - 1} (u _{T - 2}, x _{T - 1}) \\ \lambda \ge 0}} \left\{ 
            f _{T - 1} (u, x _{T - 1}) + \lambda \rho _T + \E _{\hP _T} \left[
                \sup _{x \in \cX _T} \big\{
                    Q _{T} (u, x) - \lambda c (\hx _T, x)
                \big\}
            \right]
        \right\},
    \end{align*}
    which is also bounded from below since $f _{T - 1}, \rho _T, Q _T$ are all bounded from below.
    
    Suppose we have shown $Q _{t + 1} (\cdot, \cdot)$ is random lower semi-continuous and obtain the reformulation for $Q _{t}$ that is bounded from below. Now let us show that $Q _t (\cdot, \cdot)$ is random lower semi-continuous and derive the expression for $Q _{t - 1}$ and show it is bounded from below.
    Since $Q _{t + 1} (\cdot, \cdot)$ is random lower semi-continuous, $u \mapsto Q _{t + 1} (u, x)$ is lower semi-continuous for every $x \in \cX _{t + 1}$. Therefore, for every $x,\hx_{t+1}\in\cX_{t+1}$, the function $u\mapsto Q_{t+1}(u, x) - \lambda c(\hx_{t+1},x)$ is lower semi-continuous, and thus their supremum $u\mapsto \sup_{x\in\cX_{t+1}}\left\{Q_{t+1}(u, x) - \lambda c(\hx_{t+1},x)\right\}$ is lower semi-continuous, and taking the expectation with respect to $\hP_{t+1}$, the function $\E_{\hP_{t+1}}\left[\sup_{x\in\cX_{t+1}}\big\{Q_{t+1}(u,x) - \lambda c(\hx_{t+1},x)\big\}\right]$ is lower semi-continuous in $u$. 
    It follows from \citep[Theorem 9.50]{shapiro2021lectures} that $Q_{t}(\cdot,\cdot)$ is random lower semi-continuous. Then using Example \ref{eg:suslin} we have
    \begin{align*}
        Q _{t - 1} (u _{t - 2}, x _{t - 1}) &= \inf _{u \in \cU _{t - 1} (u _{t - 2}, x _{t - 1})} \left\{
            f _{t - 1} (u, x _{t - 1}) + \sup _{\P \in \frakM _{t}} \E _{\P} \big[
                Q _{t} (u, x _{t})
            \big]
        \right\} \\
        &= \inf _{\substack{u \in \cU _{t - 1} (u _{t - 2}, x _{t - 1}) \\ \lambda \ge 0}} \left\{ 
            f _{t - 1} (u, x _{t - 1}) + \lambda \rho _t + \E _{\hP _t} \left[
                \sup _{x \in \cX_t} \big\{
                    Q _{t} (u, x) - \lambda c (\hx _t, x) 
                \big\}
            \right]
        \right\},
    \end{align*}
    which is again bounded from below, and thereby we complete the induction.
    \QED
\endproof

\proof{Proof of Example \ref{eg:01}.}
    Define $f = \mathbf1 _{\cSc}$, $c (\hx, x) = d (\hx, x) ^p$. Then \eqref{eqn:cc} is equivalent to 
    \begin{align*}
        \cL (\rho ^p) = \sup _{\P \in \cP (\cX)} \left\lbrace 
            \E _{X \sim \P} [f (X)]: \cK (\hP, \P) \le \rho ^p
        \right\rbrace \le \beta.
    \end{align*}
    For $p \in [1, \infty)$, we observe the following for each $\hx \in \cX$:
    \begin{align*}
        \sup _{x \in \cX} \{
            f (x) - \lambda c (\hx, x)
        \} 
        & = \sup _{x \in \cSc} \{
            1 - \lambda d (\hx, x) ^p
        \} \vee \sup _{x \in \cS} \{
            0 - \lambda d (\hx, x) ^p
        \} \\
        &= (1 - \lambda d (\hx, \cSc) ^p) \vee (-\lambda d (\hx, \cS) ^p) \\
        &= \begin{cases}
            1, & \hx \in \cSc \\
            (1 - \lambda d (\hx, \cSc) ^p) _+, & \hx \in \cS
        \end{cases} \\
        &= (1 - \lambda d (\hx, \cSc) ^p) _+.
    \end{align*}
    By Theorem \ref{thm}, if $f - \lambda c$ satisfies \ref{IP} for every $\lambda > 0$, then the dual problem can be calculated as the following:
    \begin{align*}
        (-\cL) ^* (-\lambda) = \E _{\hX \sim \hP} \left[ 
            \sup _{x \in \cX} \left\{
                f (x) - \lambda c (\hX, x)
            \right\}
        \right] &= \E _{\hX \sim \hP} \left[ 
            (1 - \lambda d (\hX, \cSc) ^p) _+
        \right] .
    \end{align*}
    Therefore, for every $\rho > 0$, we have 
    \begin{align*}
        \cL (\rho ^p) = \min _{\lambda \ge 0} \left\lbrace 
            \lambda \rho ^p + (-\cL) ^* (-\lambda)
        \right\rbrace = \min _{\lambda \ge 0} \left\lbrace 
            \lambda \rho ^p + \E _{\hX \sim \hP} \left[ 
                (1 - \lambda d (\hX, \cSc) ^p) _+
            \right]
        \right\rbrace.
    \end{align*}
    
    For $\beta \in (0, 1)$, the chance constraint can be written as 
    \begin{align*}
        \cL (\rho ^p) \le \beta &\iff \lambda \rho ^p + \E _{\hX \sim \hP} \left[ 
            (1 - \lambda d (\hX, \cSc) ^p) _+
        \right] \le \beta \text{ for some } \lambda \ge 0 \\
        &\iff \lambda \rho ^p + \E _{\hX \sim \hP} \left[ 
            (1 - \lambda d (\hX, \cSc) ^p) _+
        \right] \le \beta \text{ for some } \lambda > 0 \\
        &\iff \frac{\rho ^p}\beta + \frac1\beta \E _{\hX \sim \hP} \left[ 
            \left(\frac1\lambda - d (\hX, \cSc) ^p\right) _+
        \right] \le \frac1\lambda \text{ for some } \lambda > 0 \\
        &\iff \alpha + \frac1\beta \E _{\hX \sim \hP} \left[ 
            \left(-d (\hX, \cSc) ^p - \alpha\right) _+
        \right] \le -\frac{\rho ^p}\beta \text{ for some } \alpha < 0 \\
        &\iff \alpha + \frac1\beta \E _{\hX \sim \hP} \left[ 
            \left(-d (\hX, \cSc) ^p - \alpha\right) _+
        \right] \le -\frac{\rho ^p}\beta \text{ for some } \alpha \in \R \\
        &\iff \CVaRhP _\beta (-d (\hX, \cSc) ^p) = \min _{\alpha \in \R} \left\lbrace 
            \alpha + \frac1\beta \E _{\hX \sim \hP} \left[ 
                \left(-d (\hX, \cSc) ^p - \alpha\right) _+
            \right]
        \right\rbrace \le -\frac{\rho ^p}\beta.
        \tag*{\QED}
    \end{align*}
\endproof
  
\section{Proofs for Section \ref{sec:maximum-cost}}
\label{app:maximum-cost}

\proof{Proof of Proposition \ref{prop:infty-<}.}
  
We compute $\cLm ^\circ (\rho)$ as follows. 
\begin{align*}
    \cLm ^\circ (\rho) &= \sup _{\P \in \cP (\cX)} \left\lbrace
        \E _{X \sim \P} [f (X)]: \cM (\hP, \P) < \rho 
    \right\rbrace \\
    &= \sup _{\P \in \cP (\cX)} \left\lbrace
        \E _{X \sim \P} [f (X)]: \inf _{\gamma \in \Gamma (\hP, \P)} \gamma \text{-} \essssup _{\hx, x \in \cX} c (\hx, x) < \rho 
    \right\rbrace \\
    &= \sup _{\P \in \cP (\cX)} \left\lbrace
        \E _{X \sim \P} [f (X)]: \gamma \text{-} \essssup _{\hx, x \in \cX} c (\hx, x) < \rho \text{ for some } \gamma \in \Gamma (\hP, \P)
    \right\rbrace \\
    &= \sup _{\gamma \in \Gamma _{\hP}} \left\lbrace
        \E _{(\hX, X) \sim \gamma} [f (X)]: \gamma \text{-} \essssup _{\hx, x \in \cX} c (\hx, x) < \rho 
    \right\rbrace \\
    &= \sup _{\gamma \in \Gamma _{\hP}} \left\lbrace
        \E _{(\hX, X) \sim \gamma} \left[
            f (X) - \infty \mathbf1{\{c (\hX, X) \ge \rho\}}
        \right]
    \right\rbrace \\
    &= \E _{\hX \sim \hP} \left[\sup _x 
        \left\lbrace 
            f (x) - \infty \mathbf1{\{c (\hX, x) \ge \rho\}}
        \right\rbrace 
    \right].
    \\
    &= \E _{\hX \sim \hP} \left[\sup _x \left\lbrace
        f (x): c (\hX, x) < \rho
    \right\rbrace \right].
\end{align*}
In the second to the last line, we need \ref{IP} on the function $\phi (\hx, x) = f (x) - \infty \mathbf1{\{c (\hx, x) \ge \rho\}}$.

Next, we compute the dual
\begin{align*}
    (-\cLm ^\circ) ^* (-\lambda) = \sup _{\rho \ge 0} \left\lbrace
        \cLm ^\circ (\rho) - \lambda \rho 
    \right\rbrace &= \sup _{\rho \ge 0} \sup _{\P \in \cP (\cX)} \left\lbrace
        \E _{X \sim \P} [f (X)] - \lambda \rho : \cW _\infty (\hP, \P) < \rho
    \right\rbrace \\
    &= \sup _{\P \in \cP (\cX)} \left\lbrace
        \E _{X \sim \P} [f (X)] - \lambda \cW _\infty (\hP, \P) 
    \right\rbrace \\
    &= \sup _{\rho \ge 0} \sup _{\P \in \cP (\cX)} \left\lbrace
        \E _{X \sim \P} [f (X)] - \lambda \rho : \cW _\infty (\hP, \P) \le \rho
    \right\rbrace \\
    &= \sup _{\rho \ge 0} \left\lbrace
        \cLm (\rho) - \lambda \rho 
    \right\rbrace = (-\cLm) ^* (-\lambda).
\end{align*}
It follows that both $(-\cLm ^\circ) ^* (-\lambda)$ and $(-\cLm) ^* (-\lambda)$ equal the soft-constrained robust loss.
Thus
\begin{align*}
    (-\cLm) ^* (-\lambda) = (-\cLm ^\circ) ^* (-\lambda) = \sup _{\rho \ge 0} \left\lbrace 
        \E _{\hX \sim \hP} \left[\sup _x \left\lbrace
            f (x): c (\hX, x) < \rho\right] - \lambda \rho
        \right\rbrace
    \right\rbrace .
\end{align*}
This completes the proof of the proposition. 
We remark that $\cLm ^\circ$ is no longer necessarily concave, so $(-\cLm) ^{**} (\rho)$ may differ from $-\cLm ^\circ (\rho)$.
\QED 
\endproof 

\proof{Proof of Theorem \ref{thm:infty-≤}.}

Similarly to the above proof of Propositions \ref{prop:infty-<}, we have
\begin{align*}
    \cLm (\rho) &= \sup _{\P \in \cP (\cX)} \left\lbrace
        \E _{X \sim \P} [f (X)]: \cM  (\hP, \P) \le \rho 
    \right\rbrace \\
    &= \sup _{\P \in \cP (\cX)} \left\lbrace
        \E _{X \sim \P} [f (X)]: \inf _{\gamma \in \Gamma (\hP, \P)} \gamma \text{-} \essssup _{\hx, x \in \cX} c (\hx, x) \le \rho 
    \right\rbrace \\
    &\le \sup _{\P \in \cP (\cX)} \left\lbrace
        \E _{X \sim \P} [f (X)]: \gamma \text{-} \essssup _{\hx, x \in \cX} c (\hx, x) \le \rho \text{ for some } \gamma \in \Gamma (\hP, \P)
    \right\rbrace \\
    &= \sup _{\gamma \in \Gamma _{\hP}} \left\lbrace
        \E _{(\hX, X) \sim \gamma} [f (X)]: \gamma \text{-} \essssup _{\hx, x \in \cX} c (\hx, x) \le \rho 
    \right\rbrace \\
    &= \sup _{\gamma \in \Gamma _{\hP}} \left\lbrace
        \E _{(\hX, X) \sim \gamma} \left[
            f (X) - \infty \mathbf1{\{c (\hX, X) > \rho\}}
        \right]
    \right\rbrace \\
    &= \E _{\hX \sim \hP} \left[\sup _x 
        \left\lbrace 
            f (x) - \infty \mathbf1{\{c (\hX, x) > \rho\}}: c (\hX, x) \le \rho
        \right\rbrace 
    \right].
    \\
    &= \E _{\hX \sim \hP} \left[\sup _x \left\lbrace
        f (x): c (\hX, x) \le \rho
    \right\rbrace \right].
\end{align*}
Here we used \ref{IP} as it holds for all measurable functions according to Example \ref{eg:polish}.
Inequality becomes equality if 
\begin{align}
    \label{eqn:iff}
    \inf _{\gamma \in \Gamma (\hP, \P)} \gamma \text{-} \essssup _{\hx, x \in \cX} c (\hx, x) \le \rho
\end{align}
can be achieved at some $\gamma \in \Gamma (\hP, \P)$.

Now we use the additional information that $\cX$ is Polish.
If \eqref{eqn:iff} holds, we first find $\gamma _n \in \Gamma (\hP, \P)$ with $\supp \gamma _n \subset \{c \le \rho + \frac1n\}$. For any $\epsilon > 0$, we can find compact sets $\widehat K, K \subset \cX$ with $\hP [\widehat K] > 1 - \epsilon$ and $\P [K] > 1 - \epsilon$, because $\hP$ and $\P$ are probability measures on a Polish space, which are tight. Then $\gamma _n [\widehat K \times K] > 1 - 2 \epsilon$ for each $n$. This shows $\{\gamma _n\} _n$ is a tight sequence. Since $\cX$ is complete and separable, by Prokhorov theorem, there exists a weakly converging subsequence $\gamma _{n _k} \to \gamma \in \cP (\cX \times \cX)$. Marginals of $\gamma _n$ also converge weakly to the marginals of $\gamma$, so $\gamma \in \Gamma (\hP, \P)$. To show that $\gamma$ is supported in $\{c \le \rho\}$, define $g: \cX \times \cX \to \R$ by  
\begin{align*}
    g (\hx, x) = 1 \wedge (c (\hx, x) - \rho) _+.
\end{align*}
$g$ is continuous and bounded on $\cX \times \cX$. Moreover, $\E _{\gamma _n} [g] \le \frac1n$. By weak convergence, $\E _{\gamma} [g] = 0$, so $c (\hx, x) \le \rho$ for $\gamma$-a.e. $(\hx, x) \in \cX \times \cX$. That is, $\gamma\text{-}\essssup c \le \rho$.

Next, we compute the dual.
\begin{align*}
    (-\cLm) ^* (-\lambda) = \sup _{\rho \ge 0} \left\lbrace
        \cLm (\rho) - \lambda \rho 
    \right\rbrace &= \sup _{\rho \ge 0} \sup _{\P \in \cP (\cX)} \left\lbrace
        \E _{X \sim \P} [f (X)] - \lambda \rho : \cW _\infty (\hP, \P) \le \rho
    \right\rbrace \\
    &= \sup _{\P \in \cP (\cX)} \left\lbrace
        \E _{X \sim \P} [f (X)] - \lambda \cW _\infty (\hP, \P) 
    \right\rbrace.
\end{align*}
Thus
\begin{align*}
    (-\cLm) ^* (-\lambda) = \sup _{\rho \ge 0} \left\lbrace 
        \E _{\hX \sim \hP} \left[\sup _x \left\lbrace
            f (x): c (\hX, x) \le \rho\right] - \lambda \rho
        \right\rbrace
    \right\rbrace .
\end{align*}
This completes the proof of the theorem.
\QED
\endproof

\proof{Proof of Example \ref{eg:01-infinity}.}
    By Theorem \ref{thm:infty-≤}, we have 
    \begin{align*}
        \cLm (\rho) = \sup _{\P \in \cP (\cX)} \{
            \E _{X \sim \P} [f (X)]: \cW _\infty (\hP, \P) \le \rho
        \} = \E _{\hX \sim \hP} \left[
            \sup _x \left\lbrace
            f (x): d(\hX, x) \le \rho
        \right\rbrace \right] = \hP (d (\hX, \cSc) \le \rho).
    \end{align*}
    In particular, $\cLm (\rho) = 1$ if $\rho \ge d (\supp \hP, \cSc)$.
    We remark that now the corresponding soft robust problem is
    \begin{align*}
        (-\cLm) ^* (-\lambda) = \sup _{\rho \ge 0} \left\lbrace 
            \hP (d (\hX, \cSc) \le \rho) - \lambda \rho
        \right\rbrace .
        \tag*{\QED}
    \end{align*}
\endproof

\section{Proofs for Section \ref{sec:risk}}
\label{app:risk}

\proof{Proof of Theorem \ref{thm:risk-measure}.}

First, we consider the maximum transport cost $\cM$. 
Denote 
\begin{align*}
    J' (\P) := \inf _{\alpha \in A'} \E _{X \sim \P} [f _\alpha (X)].
\end{align*}
By assumption \eqref{eqn:achievable}, we know that 
\begin{align*}
    \cLJm (\rho) = \sup _{\P \in \cP (\cX)} \left\lbrace
        J' (\P): \cM (\hP, \P) \le \rho
    \right\rbrace
\end{align*}
Similar as Theorem \ref{thm:infty-≤}, we have 
\begin{align*}
    \cLJm (\rho) 
    &= \sup _{\P \in \cP (\cX)} \left\lbrace
        \inf _{\alpha \in A'} \E _{X \sim \P} [f _\alpha (X)]: \cM (\hP, \P) \le \rho 
    \right\rbrace \\
    &= \sup _{\P \in \cP (\cX)} \left\lbrace
        \inf _{\alpha \in A'} \E _{X \sim \P} [f _\alpha (X)]: \inf _{\gamma \in \Gamma (\hP, \P)} \gamma \text{-} \essssup _{\hx, x \in \cX} c (\hx, x) \le \rho 
    \right\rbrace \\
    &= \sup _{\P \in \cP (\cX)} \left\lbrace
        \inf _{\alpha \in A'} \E _{X \sim \P} [f _\alpha (X)]: \gamma \text{-} \essssup _{\hx, x \in \cX} c (\hx, x) \le \rho \text{ for some } \gamma \in \Gamma (\hP, \P)
    \right\rbrace \\
    &= \sup _{\gamma \in \Gamma _{\hP}} \left\lbrace
        \inf _{\alpha \in A'} \E _{X \sim \P} [f _\alpha (X)]: \gamma \text{-} \essssup _{\hx, x \in \cX} c (\hx, x) \le \rho 
    \right\rbrace \\
    &= \sup _{\gamma \in \Gamma _{\hP}} \left\lbrace
        \inf _{\alpha \in A'} \E _{(\hX, X) \sim \gamma} \left[
            f _\alpha (X) - \infty \mathbf1{\{c (\hX, X) > \rho\}}
        \right]
    \right\rbrace.
\end{align*}
We claim that for each $\gamma \in \Gamma _{\hP}$, $A' \ni \alpha \mapsto \E _{\gamma} [f _\alpha (X)]$ has the following properties:
\begin{enumerate}
    \item[(a)] Convexity: given $\alpha _0, \alpha _1 \in \R$, $\theta \in (0, 1)$, define $\alpha _\theta = (1 - \theta) \alpha _0 + \theta \alpha _1$. Due to convexity of $f _\alpha$ in $\alpha$ we have $f _{\alpha _\theta} (x) \le (1 - \theta) f _{\alpha _0} (x) + \theta f _{\alpha _1} (x)$ for every $x \in \cX$. So 
    \begin{align*}
        \E _{\gamma} [f _{\alpha _\theta} (X)] \le (1 - \theta) \E _{\gamma} [f _{\alpha _0} (X)] + \theta \E _{\gamma} [f _{\alpha _1} (X)].
    \end{align*}

    \item[(b)] Lower semi-continuity: let $\alpha _n \to \alpha$ as $n \to \infty$. Due to lower semi-continuity of $f _\alpha$ in $\alpha$ we have $f _\alpha (x) \le \liminf _{n \to \infty} f _{\alpha _n} (x)$ for every $x \in \cX$. So 
    \begin{align*}
        \E _{\gamma} [f _{\alpha} (X)] \le \E _{\gamma} \left[\liminf _{n \to \infty} f _{\alpha _n} (X)\right] \le \liminf _{n \to \infty} \E _{\gamma} [f _{\alpha _n} (X)].
    \end{align*}
    The second inequality is due to Fatou's lemma thanks to $\{f _\alpha\} _{\alpha \in A'}$ being uniformly bounded from below.
\end{enumerate}
Since $A'$ is compact, by Sion's minimax theorem and interchangeability, 
\begin{align*}
    \cLJm (\rho) 
    &= \inf _{\alpha \in A'}  \left\lbrace
        \sup _{\gamma \in \Gamma _{\hP}} \E _{(\hX, X) \sim \gamma} \left[
            f _\alpha (X) - \infty \mathbf1{\{c (\hX, X) > \rho\}}
        \right]
    \right\rbrace \\
    &= \inf _{\alpha \in A'} \E _{\hX \sim \hP} \left[\sup _x 
        \left\lbrace 
            f _\alpha (x) - \infty \mathbf1{\{c (\hX, x) > \rho\}}: c (\hX, x) \le \rho
        \right\rbrace 
    \right].
    \\
    &= \inf _{\alpha \in A'} \E _{\hX \sim \hP} \left[\sup _x \left\lbrace
        f _\alpha (x): c (\hX, x) \le \rho
    \right\rbrace \right] =: \inf _{\alpha \in A'} \ell _\alpha.
\end{align*}
To complete the proof, we need to enlarge $A'$ to $A$ again. Let $\alpha _n$ be a minimizing sequence:
\begin{align*}
    \lim _{n \to \infty} \ell _{\alpha _n} = \inf _{\alpha \in A} \ell _\alpha.
\end{align*}
Denote $A' _n = \operatorname{conv} (A' \cup \{\alpha _n\})$ to be the convex hull of $A'$ and $\alpha _n$. Clearly, if \eqref{eqn:achievable} holds for $A'$, it should also hold for $A' _n \supset A'$. Since $A' _n$ is still compact, the same argument on $A' _n$ instead of $A'$ gives $\cLJm (\rho) = \inf _{\alpha \in A' _n} \ell _\alpha$. This holds for every $n$, so taking the limit yields
\begin{align*}
    \cLJm (\rho) = \inf _{\alpha \in A} \ell _\alpha &= \inf _{\alpha \in A} \E _{\hX \sim \hP} \left[\sup _x \left\lbrace
        f _\alpha (x): c (\hX, x) \le \rho
    \right\rbrace \right].
\end{align*}

Next, we consider the Kantorovich cost $\cK$. Similarly, we can restrict to $A'$ by assumption \eqref{eqn:achievable}:
\begin{align*}
    \cLJ (\rho) = \sup _{\P \in \cP (\cX)} \left\lbrace
        J' (\P): \cK (\hP, \P) \le \rho
    \right\rbrace.
\end{align*}
It is easy to see that $J'$ is concave in $\P$.
Similar as Lemma \ref{lem:concave}, it is a simple exercise to show $\cLJ (\cdot)$ is lower bounded by $\sup _{\alpha \in A'} \E_{\hP} [f _\alpha]$, monotonically increasing, and concave on $[0,\infty)$. Now we take the dual of $-\cLJ$: 
\begin{align*}
    (-\cLJ) ^* (-\lambda) := \sup _{\rho \ge 0} \left\lbrace 
        \cLJ (\rho) - \lambda \rho
    \right\rbrace &= \sup _{\rho \ge 0, \P \in \cP (\cX)} \left\lbrace 
        J' (\P) - \lambda \rho: \cK (\hP, \P) \le \rho
    \right\rbrace \\
    &= \sup _{\P \in \cP (\cX)} \left\lbrace 
        J' (\P) - \lambda \cK (\hP, \P)%
    \right\rbrace \\
    &= \sup _{\gamma \in \Gamma _{\hP}} \inf _{\alpha \in A'} \left\lbrace 
        \E _{\gamma} [f _\alpha (X) - \lambda c (\hX, X)]
    \right\rbrace.
\end{align*}

We have shown that $\alpha \mapsto \E _{\gamma} [f _\alpha (X)]$ is lower semi-continuous and convex in the first half of the proof. By Sion's minimax theorem, we can exchange sup and inf, so 
\begin{align*}
    (-\cLJ) ^* (-\lambda)
    &= \inf _{\alpha \in A'} \sup _{\gamma \in \Gamma _{\hP}} \left\lbrace 
        \E _{\gamma} [f _\alpha (X) - \lambda c (\hX, X)]
    \right\rbrace = \inf _{\alpha \in A'} \E _{\hX \sim \hP} \left[ 
        \sup _{x \in \cX} \left\lbrace 
            f _\alpha (x) - \lambda c (\hX, x)
        \right\rbrace
    \right].
\end{align*}
Here we used \ref{IP} on $\phi _{\lambda, \alpha} = f _\alpha - \lambda c$. By enlarging $A'$ to $A$ as in the maximal cost case, we have
\begin{align*}
    (-\cLJ) ^* (-\lambda)  = \inf _{\alpha \in A} \E _{\hX \sim \hP} \left[ 
        \sup _{x \in \cX} \left\lbrace 
            f _\alpha (x) - \lambda c (\hX, x)
        \right\rbrace
    \right].
    \tag*{\QED}
\end{align*}

\endproof

To apply Theorem \ref{thm:risk-measure} on Example \ref{eg:risk}, we need to verify that the infimum is achieved in a finite interval dependent only on the transport cost. The following lemma confirms this property for $\CVaR$ and $\MAD$.

\begin{lemma}
    \label{lem:VaR}
    Suppose $\cX = \R$.
    Let $\hP, \P \in \cP (\cX)$, $\beta \in (0, 1)$. Let $\hX \sim \hP$, $X \sim \P$. If $\P (X \ge \alpha) \ge \beta$ and $\P (X \le \alpha) \ge 1 - \beta$, then 
    \begin{align*}
         \alpha \in \left[
            -\CVaRhP _{1 - \beta} (-\hX) - \frac1{1 - \beta} \cW _1 (\hP, \P)    
            ,
            \CVaRhP _\beta (\hX) + \frac1\beta \cW _1 (\hP, \P)
        \right].
    \end{align*}
\end{lemma}

\proof{Proof.}

Given any $\gamma \in \Gamma (\hP, \P)$, we have 
\begin{align*}
    \E _{(\hX, X) \sim \gamma} \left[\|\hX - X\|\right] &\ge \E _{(\hX, X) \sim \gamma} \left[\|X - \hX\| \mathbf1\{X \ge \alpha\}\right] \\
    &\ge \E _{(\hX, X) \sim \gamma} \left[(\alpha - \hX) \mathbf1\{X \ge \alpha\}\right] \\
    &= \alpha \P (X \ge \alpha) - \E _{\hX \sim \hP} \left[ 
        \hX \E _{X \sim \gamma _{X \mid \hX}} \left[ 
            \mathbf1\{X \ge \alpha\} \mid \hX
        \right]
    \right] \\
    &= \P (X \ge \alpha) \left(
        \alpha - \E _{\hX \sim \hP} \left[ 
            \hX \cdot \frac{\P (X \ge \alpha \mid \hX)}{\P (X \ge \alpha)}
        \right]
    \right)
\end{align*}
Note that $\E _{\hX \sim \hP} \left[\dfrac{\P (X \ge \alpha \mid \hX)}{\P (X \ge \alpha)}\right] = 1$, and $\dfrac{\P (X \ge \alpha \mid \hX)}{\P (X \ge \alpha)} \le \dfrac{1}{\beta}$. Therefore 
\begin{align*}
    \E _{\hX \sim \hP} \left[ 
        \hX \cdot \frac{\P (X \ge \alpha \mid \hX)}{\P (X \ge \alpha)}
    \right] \le \sup _{\hQ \in \cP (\cX)} \left\lbrace 
        \E _{\hX \sim \hQ} [\hX] : \hQ \ll \hP, \frac{d \hQ}{d \hP} \le \frac1\beta
    \right\rbrace = \CVaRhP _\beta (\hX).
\end{align*}
Here we used the dual formulation for CVaR in \cite{follmer2010convex}. Hence, we have shown that
\begin{align*}
    \alpha - \CVaRhP _\beta (\hX) \le \frac1{\P (X \ge \alpha)} \E _{(\hX, X) \sim \gamma} \left[\|\hX - X\|\right] \le \frac1\beta \E _{(\hX, X) \sim \gamma} \left[\|\hX - X\|\right].
\end{align*}
The proof of the upper bound is completed by taking infimum over $\gamma \in \Gamma (\hP, \P)$. The proof of the lower bound is similar:
\begin{align*}
    \E _{(\hX, X) \sim \gamma} \left[\|\hX - X\|\right] 
    &\ge \E _{(\hX, X) \sim \gamma} \left[\|\hX - X\| \mathbf1\{X \le \alpha\}\right] \\
    &\ge \E _{(\hX, X) \sim \gamma} \left[(\hX - \alpha) \mathbf1\{X \le \alpha\}\right] \\
    &= \E _{\hX \sim \hP} \left[ 
        \hX \E _{X \sim \gamma _{X \mid \hX}} \left[ 
            \mathbf1\{X \le \alpha\} \mid \hX
        \right]
    \right] - \alpha \P (X \le \alpha) \\
    &= \P (X \le \alpha) \left(
        \E _{\hX \sim \hP} \left[ 
            \hX \cdot \frac{\P (X \le \alpha \mid \hX)}{\P (X \le \alpha)}
        \right] - \alpha
    \right) \\
    &\ge \P (X \le \alpha) \left(
        -\CVaRhP _{1 - \beta} (-\hX) - \alpha
    \right). 
    \tag*{\QED}
\end{align*}
\endproof 

\proof{Proof of Example \ref{eg:risk} ($\CVaR$).}%

Given $\hZ \sim \hQ$ and $Z \sim \Q$, define $\hX = b ^\top \hZ$ and $X = b ^\top Z$, and let $\hP$, $\P$ denote the law of $\hX$ and $X$ respectively. We observe the following 
\begin{enumerate}
    \item[(a)] $\CVaRQ _\beta (b ^\top Z) = \CVaRP _\beta (X)$ and $\CVaRhQ _\beta (b ^\top \hZ) = \CVaRhP _\beta (\hX)$.

    \item[(b)] For any $\Q \in \cP (\cZ)$, $\cW _p (\hP, \P) \le \| b \| _* \cW _p (\hQ, \Q)$.

    \item[(c)] For any $\P' \in \cP (\R)$, we can find $\Q \in \cP (\cZ)$ such that $\P = \P'$, and $\cW _p (\hP, \P) \ge \| b \| _* \cW _p (\hQ, \Q)$.
\end{enumerate}
If all these claims are true, then 
\begin{align*}
    \sup _{\Q \in \cP (\cZ)} \left\{
        \CVaRQ _\beta (b ^\top Z): \cW _p (\hQ, \Q) \le \rho
    \right\} = \sup _{\P \in \cP (\R)} \left\{
        \CVaRP _\beta (X): \cW _p (\hP, \P) \le \| b \| _* \rho
    \right\}.
\end{align*}
The first two claims are direct: $X = b ^\top Z$, $\hX = b ^\top \hZ$, and $|X - \hX| \le \| b \| _* \| Z - \hZ \|$. For the third claim, we prove it as follows. Let $b ^* \in \cZ$ be the unit dual of $b ^\top$, i.e. $b ^\top b ^* = \| b ^\top \| _*$ and $\|b ^*\| = 1$. Given $\hZ \sim \hP$, $\hX = b ^\top X \sim \hQ$, $X' \sim \P'$, define $Z = \hZ + (X' - \hX) b ^* / \| b ^\top \| _*$. Then $X = b ^\top Z = X$, $\P' = \P$, and $\| \hZ - Z \| = \| b ^\top \| _* ^{-1} \| \hX - X \|$, thus $\cW _p (\hQ, \Q) \le \| b ^\top \| _* ^{-1} \cW _p (\hP, \P)$. 
We have transformed the problem to the following form: with $p < \infty$, $c (\hx, x) = |\hx - x| ^p$, 
\begin{align*}
    \sup _{\Q \in \cP (\cZ)} \left\lbrace
        \CVaRQ _\beta (b ^\top Z): \cW _p (\hQ, \Q) \le \rho
    \right\rbrace = 
    \sup _{\P \in \cP (\cX)} \left\{
        J (\P): \cK (\hP, \P) \le (\| b ^\top \| _* \rho) ^p
    \right\} = \cLJ ((\| b ^\top \| _* \rho) ^p),
\end{align*} 
and with $p = \infty$, $c (\hx, x) = |\hx - x| ^p$,  
\begin{align*}
    \sup _{\Q \in \cP (\cZ)} \left\lbrace
        \CVaRQ _\beta (b ^\top Z): \cW _p (\hQ, \Q) \le \rho
    \right\rbrace = 
    \sup _{\P \in \cP (\cX)} \left\{
        J (\P): \cM (\hP, \P) \le \| b ^\top \| _* \rho
    \right\} = \cLJm (\| b ^\top \| _* \rho).
\end{align*}

For simplicity, assume $\| b ^* \| = 1$ from now on.

To apply Theorem \ref{thm:risk-measure}, we verify the following prerequisites:

\begin{itemize}
    \item $f _\alpha$ satisfies Assumption \ref{assum:setup}: $f _\alpha \ge \alpha$ so $\E _{\hP} [f _\alpha] \ge \alpha > -\infty$.
    \item $f _\alpha$ is lower semi-continuous and convex in $\alpha$: this is obvious since $f _\alpha = \max \{\alpha, \frac1\beta x + (1 - \frac1\beta) \alpha \}$ is the maximum of two affine functions.
    \item $\inf _{\alpha \in A', x \in \cX} f _\alpha (x) > -\infty$ for compact $A' \subset \R$: $\inf _{\alpha \in A', x \in \cX} f _\alpha (x) = \inf _{\alpha \in A'} \alpha = \min _{\alpha \in A'} \alpha > -\infty$, since $A'$ is compact and bounded.
    \item $f _\alpha - \lambda c$ satisfies \ref{IP}: this is because $\cX$ is a Euclidean space, which is complete and separable. 
    \item \eqref{eqn:achievable} holds for $\P$ in Wasserstein ball: from \cite{rockafellar2000optimization} we know that for CVaR problem, the minimum of \eqref{eqn:JP} is attained on a nonempty closed bounded interval $\alpha \in [\underline \alpha, \overline \alpha]$ (possibly a singleton). This interval contains $\alpha$ such that $\P (X \ge \alpha) \ge \beta$ and $\P (X \le \alpha) \ge 1 - \beta$. By Lemma \ref{lem:VaR}, this interval is contained in 
    \begin{align*}
        A' = \left[
            -\CVaRhP _{1 - \beta} (-\hX) - \frac1{1 - \beta} \cW _1 (\hP, \P)    
            ,
            \CVaRhP _\beta (\hX) + \frac1\beta \cW _1 (\hP, \P)
        \right].
    \end{align*}
\end{itemize}
Since $\cW _1 (\hP, \P) \le \cW _p (\hP, \P)$ for $p \in [1, \infty]$, we have verified all the prerequisites of Theorem \ref{thm:risk-measure}.

When $p = \infty$, by Theorem \ref{thm:risk-measure}, we conclude
\begin{align*}
    \cLJm (\rho) = \sup _{\P \in \cP} \left\lbrace
        \CVaRP _\beta (X): \cW _\infty (\hP, \P) \le \rho
    \right\rbrace = \inf _{\alpha \in A} \E _{\hX \sim \hP} \left[\sup _x \left\lbrace
        f _\alpha (x): |\hX - x| \le \rho
    \right\rbrace \right],
\end{align*}
where
\begin{align*}
    \sup _x \left\lbrace
        f _\alpha (x): |\hX - x| \le \rho
    \right\rbrace = \sup _{x} \left\lbrace
        \alpha + \frac1{1 - \beta} (x - \alpha) _+: |\hX - x| \le \rho
    \right\rbrace = \alpha + \frac1{1 - \beta} (\hX + \rho - \alpha) _+.
\end{align*}
We thus conclude that  
\begin{align*}
    \cLJm (\rho) = \inf _{\alpha \in A} \left\{
        \alpha + \frac1{1 - \beta} \E _{\hX \sim \hP} [(\hX + \rho - \alpha) _+]
    \right\} = \CVaRhP _\beta (\hX) + \rho.
\end{align*}

When $p = 1$, 
\begin{align*}
    \sup _{x \in \cX} \left\lbrace 
        f _\alpha (x) - \lambda c (\hx, x)
    \right\rbrace = \sup _{x \in \cX} \left\lbrace 
        \alpha + \frac1{1 - \beta} (x - \alpha) _+ - \lambda |\hx - x|
    \right\rbrace = \begin{cases}
        \displaystyle
        \alpha + \frac1{1 - \beta} (\hx - \alpha) _+ & \lambda \ge \frac1{1 - \beta} \\
        \infty & \lambda < \frac1{1 - \beta}
    \end{cases}.
\end{align*}
Therefore for $\lambda \ge \frac1{1 - \beta}$,
\begin{align*}
    \inf _{\alpha \in A} \left\lbrace 
        \E _{\hX \sim \hP} \left[ 
            \sup _{x \in \cX} \left\lbrace 
                f _\alpha (x) - \lambda c (\hX, x)
            \right\rbrace
        \right]
    \right\rbrace = \inf _{\alpha \in A} \left\lbrace 
        \E _{\hX \sim \hP} \left[ 
            \alpha + \frac1{1 - \beta} (\hX - \alpha) _+
        \right]
    \right\rbrace = \CVaRhP _\beta (\hX).
\end{align*}
Thus
\begin{align*}
    \cLJ (\rho) = \inf _{\alpha \in A, \lambda \ge 0} \left\lbrace 
        \lambda \rho + \E _{\hX \sim \hP} \left[ 
            \sup _{x \in \cX} \left\lbrace 
                f _\alpha (x) - \lambda c (\hX, x)
            \right\rbrace
        \right]
    \right\rbrace = \CVaRhP _\beta (\hX) + \frac \rho{1 - \beta}.
\end{align*}

When $p > 1$, 
\begin{align*}
    \sup _{x \in \cX} \left\lbrace 
        f _\alpha (x) - \lambda |\hx - x| ^p
    \right\rbrace &= \sup _{x \in \cX} \left\lbrace 
        \alpha + \frac1{1 - \beta} (x - \alpha) _+ - \lambda |\hx - x| ^p
    \right\rbrace \\
    &= \sup _{x \in \cX} \left\lbrace 
        \alpha + \frac1{1 - \beta} (x - \alpha) - \lambda |\hx - x| ^p
    \right\rbrace \vee \sup _{x \in \cX} \left\lbrace 
        \alpha - \lambda |\hx - x| ^p
    \right\rbrace \\
    &= \left( 
        \alpha + \frac1{1 - \beta} (\hx - \alpha) + \sup _{t \in \R} \left\lbrace 
            \frac t{1 - \beta} - \lambda |t| ^p
        \right\rbrace 
    \right) \vee \alpha \\
    &= \alpha + \left(\frac1{1 - \beta} (\hx - \alpha) + C \lambda ^{-\frac1{p - 1}}\right) _+ \\
    &= \alpha + \left(
        \frac1{1 - \beta} \left(\hx - (\alpha - C (1 - \beta) \lambda ^{-\frac1{p - 1}}\right)
    \right) _+,
\end{align*}
where $C = (p - 1) (p (1 - \beta)) ^{-\frac p{p - 1}}$. Thus 
\begin{align*}
    \inf _{\alpha \in A} \left\lbrace 
        \E _{\hX \sim \hP} \left[ 
            \sup _{x \in \cX} \left\lbrace 
                f _\alpha (x) - \lambda c (\hX, x)
            \right\rbrace
        \right]
    \right\rbrace = \CVaRhP _\beta (\hX) + C (1 - \beta) \lambda ^{-\frac1{p - 1}}. 
\end{align*}
Therefore
\begin{align*}
    \cLJ (\rho ^p) = \CVaRhP _\beta (\hX) + \min _{\lambda \ge 0} \left\lbrace 
        \lambda \rho ^p + C (1 - \beta) \lambda ^{-\frac1{p - 1}}
    \right\rbrace = \CVaRhP _\beta (\hX) + \rho (1 - \beta) ^{-\frac1p}.
\end{align*}

In conclusion, for $p \in [1, \infty]$, it holds that 
\begin{align}
    \sup _{\P \in \cP} \left\lbrace
        \CVaRP _\beta (X): \cW _p (\hP, \P) \le \rho
    \right\rbrace = \CVaRhP _\beta (\hX) + \rho (1 - \beta) ^{-\frac1p}.
    \tag*{\QED}
\end{align}
\endproof

\proof{Proof of Example \ref{eg:risk} ($\Var$).}%
    Same as Example \ref{eg:risk} $\CVaR$, we can reduce to a one-dimensional problem and assume without loss of generality that $\| b ^\top \| _* = 1$.
    
    It is well-known that the optimal $\alpha$ is the expectation:
    \begin{align*}
        \VarP (X) = \min _\alpha \E [(X - \alpha) ^2] = \E [(X - \E [X]) ^2].
    \end{align*}
    Given $\hP, \P$ and a transport plan $\gamma \in \Gamma (\hP, \P)$, the transport cost
    \begin{align*}
        \E _\gamma [\|\hX - X\|] \ge |\E _\gamma [\hX - X]| \ge |\E _{\hP} [\hX] - \E _{\P} [X]|.
    \end{align*}
    Minimizing over all $\gamma \in \Gamma (\hP, \P)$ gives 
    \begin{align*}
        \E _{X \sim \P} [X] \le \left[
            \E _{\hX \sim \hP} [\hX] - \cW _1 (\hX, X), \E _{\hX \sim \hP} [\hX] + \cW _1 (\hX, X)
        \right].
    \end{align*}
    Same as before, we can verify that $f _\alpha (x) = (x - \alpha) ^2$ meets all the prerequisites of Theorem \ref{thm:risk-measure}.
    
    When $p = \infty$, by Theorem \ref{thm:risk-measure}, we conclude
    \begin{align*}
        \cLJm (\rho) = \sup _{\P \in \cP} \left\lbrace
            \VarP (X): \cW _\infty (\hP, \P) \le \rho
        \right\rbrace = \inf _{\alpha \in A} \E _{\hX \sim \hP} \left[\sup _x \left\lbrace
            f _\alpha (x): |\hX - x| \le \rho
        \right\rbrace \right],
    \end{align*}
    where
    \begin{align*}
        \sup _x \left\lbrace
            f _\alpha (x): |\hX - x| \le \rho
        \right\rbrace = \sup _{x} \left\lbrace
            (x - \alpha) ^2: |\hX - x| \le \rho
        \right\rbrace = (|\hX - \alpha| + \rho) ^2.
    \end{align*}
    We thus conclude that  
    \begin{align*}
        \cLJm (\rho) = \inf _{\alpha \in A} \left\{
            \E _{\hX \sim \hP} \left[
                (|\hX - \alpha| + \rho) ^2
            \right]
        \right\}.
    \end{align*}
    
    When $1 \le p < 2$, for any $\lambda \ge 0$, it holds that 
    \begin{align*}
        \sup _{x \in \cX} \left\lbrace 
            f _\alpha (x) - \lambda c (\hx, x)
        \right\rbrace = \sup _{x \in \cX} \left\lbrace 
            (x - \alpha) ^2 - \lambda |\hx - x| ^p
        \right\rbrace = +\infty.
    \end{align*}
    Thus for any $\rho > 0$ we must have 
    \begin{align*}
        \cLJ (\rho ^p) = +\infty.
    \end{align*}
    
    When $p = 2$, 
    \begin{align*}
        \sup _{x \in \cX} \left\lbrace 
            f _\alpha (x) - \lambda |\hx - x| ^2
        \right\rbrace &= \sup _{x \in \cX} \left\lbrace 
            (x - \alpha) ^2 - \lambda (\hx - x) ^2
        \right\rbrace \\
        &= \begin{cases}
            + \infty & 0 \le \lambda < 1, \text{ or } \lambda = 1, \hx \neq \alpha, \text{ or } \\
            0 & \lambda = 1, \hx = \alpha \\
            \frac{\lambda}{\lambda - 1} (\hx - \alpha) ^2 & \lambda > 1.
        \end{cases}
    \end{align*}
    Thus 
    \begin{align*}
        \inf _{\alpha \in A} \left\lbrace 
            \E _{\hX \sim \hP} \left[ 
                \sup _{x \in \cX} \left\lbrace 
                    f _\alpha (x) - \lambda c (\hX, x)
                \right\rbrace
            \right]
        \right\rbrace = \begin{cases}
            + \infty & 0 \le \lambda < 1, \text{ or } \lambda = 1, \hP \neq \Dirac _{\hx} \text{ for any } \hx \in \cX \\
            0 & \lambda = 1, \hP = \Dirac _{\hx} \text{ for some } \hx \in \cX \\
            \frac{\lambda}{\lambda - 1} \VarhP (\hX) & \lambda > 1.
        \end{cases}
    \end{align*}
    We now conclude 
    \begin{align*}
        \cLJ (\rho ^2) = \inf _{\lambda > 1} \left\{
            \lambda \rho ^2 + \frac{\lambda}{\lambda - 1} \VarhP (\hX) 
        \right\} = (\VarhP (\hX) ^\frac12 + \rho) ^2. 
        \tag*{\QED}
    \end{align*}
\endproof

\proof{Proof of Example \ref{eg:risk} ($\MAD$).}%
    Same as Example \ref{eg:risk} $\CVaR$, we can reduce to a one-dimensional problem and assume without loss of generality that $\| b ^\top \| _* = 1$.
    
    It is well-known that the optimal $\alpha$ is the median. By Lemma \ref{lem:VaR}, the median of $\P$ inside the Wasserstein uncertainty set is attained in 
    \begin{align*}
         \alpha \in \left[
            -\CVaRhP _{\frac12} (-\hX) - 2 \cW _1 (\hP, \P)    
            ,
            \CVaRhP _\frac12 (\hX) + 2 \cW _1 (\hP, \P)
        \right].
    \end{align*} 
    Same as in Example \ref{eg:risk} $\CVaR$, we can verify that $f _\alpha (x) = |x - \alpha|$ meets all the prerequisites of Theorem \ref{thm:risk-measure}.

    When $p = \infty$, by Theorem \ref{thm:risk-measure}, we conclude
    \begin{align*}
        \cLJm (\rho) = \sup _{\P \in \cP} \left\lbrace
            \MADP (X): \cW _\infty (\hP, \P) \le \rho
        \right\rbrace = \inf _{\alpha \in A} \E _{\hX \sim \hP} \left[\sup _x \left\lbrace
            f _\alpha (x): |\hX - x| \le \rho
        \right\rbrace \right],
    \end{align*}
    where
    \begin{align*}
        \sup _x \left\lbrace
            f _\alpha (x): |\hX - x| \le \rho
        \right\rbrace = \sup _{x} \left\lbrace
            |x - \alpha|: |\hX - x| \le \rho
        \right\rbrace = |\hX - \alpha| + \rho.
    \end{align*}
    We thus conclude that  
    \begin{align*}
        \cLJm (\rho) = \inf _{\alpha \in A} \left\{
            \E _{\hX \sim \hP} [|\hX - \alpha|] + \rho
        \right\} = \MADhP (\hX) + \rho.
    \end{align*}
    
    When $p = 1$, 
    \begin{align*}
        \sup _{x \in \cX} \left\lbrace 
            f _\alpha (x) - \lambda c (\hx, x)
        \right\rbrace = \sup _{x \in \cX} \left\lbrace 
            |x - \alpha| - \lambda |\hx - x|
        \right\rbrace = \begin{cases}
            \displaystyle
            |\hx - \alpha| & \lambda \ge 1 \\
            \infty & \lambda < 1
        \end{cases}.
    \end{align*}
    Therefore for $\lambda \ge 1$,
    \begin{align*}
        \inf _{\alpha \in A} \left\lbrace 
            \E _{\hX \sim \hP} \left[ 
                \sup _{x \in \cX} \left\lbrace 
                    f _\alpha (x) - \lambda c (\hX, x)
                \right\rbrace
            \right]
        \right\rbrace = \inf _{\alpha \in A} \left\lbrace 
            \E _{\hX \sim \hP} \left[ 
                |\hX - \alpha|
            \right]
        \right\rbrace = \MADhP (\hX).
    \end{align*}
    Thus
    \begin{align*}
        \cLJ (\rho) = \inf _{\alpha \in A, \lambda \ge 0} \left\lbrace 
            \lambda \rho + \E _{\hX \sim \hP} \left[ 
                \sup _{x \in \cX} \left\lbrace 
                    f _\alpha (x) - \lambda c (\hX, x)
                \right\rbrace
            \right]
        \right\rbrace = \MADhP (\hX) + \rho.
    \end{align*}

    Robust loss for $p = 1$ and $p = \infty$ are both $\MADhP (\hX) + \rho$. As $\cW _1 (\hP, \P) \le \cW _p (\hP, \P) \le \cW _\infty (\hP, \P)$, we have for any $1 \le p \le \infty$:
    \begin{align}
        \sup _{\P \in \cP} \left\lbrace
            \CVaRP _\beta (X): \cW _p (\hP, \P) \le \rho
        \right\rbrace = \MADhP (\hX) + \rho.
        \tag*{\QED}
    \end{align}
\endproof

\proof{Proof of Example \ref{eg:risk} ($\Ent$).}%
    Same as Example \ref{eg:risk} $\CVaR$, we can reduce to a one-dimensional problem and assume without loss of generality that $\| b ^\top \| _* = 1$.
    
    $\alpha \mapsto \E [f _\alpha (X)] = \alpha + \frac1\theta \left(\E [e ^{\theta (X - \alpha)}] - 1\right)$ is convex, and $\lim _{\alpha \to \pm \infty} \E [f _\alpha (X)] = +\infty$, so the minimizer $\alpha ^*$ satisfies 
    \begin{align*}
        0 = \frac{d}{d \alpha} \bigg| _{\alpha = \alpha ^*} \E [f _\alpha (X)] = 1 - \E [e ^{\theta (X - \alpha ^*)}].%
    \end{align*}
    Therefore, $e ^{\theta \alpha ^*} = \E [e ^{\theta X}]$. Suppose $\halpha$ is the minimizer to $\E [f _{\alpha} (\hX)]$, and $\gamma\text{-}\essssup _{\hx, x} \| \hx - x \| \le \rho$, then 
    \begin{align*}
        e ^{\theta \alpha ^*} = \E [e ^{\theta X}] \le \E [e ^{\theta (\hX + \rho)}] = \E [e ^{\theta \hX}] e ^{\theta \rho} = e ^{\theta (\halpha + \rho)}.
    \end{align*}
    Hence $\alpha ^* \le \halpha + \rho$. Similarly, $\alpha ^* \ge \halpha - \rho$. Therefore, 
    \begin{align*}
        \alpha \in [\halpha - \cW _{\infty} (\hP, \P), \halpha + \cW _{\infty} (\hP, \P)].
    \end{align*}
    By Theorem \ref{thm:risk-measure}, we conclude for $p = \infty$:
    \begin{align*}
        \cLJm (\rho) &= \inf _{\alpha \in \R} \E _{\hX \sim \hP} \left[ 
            \sup _{x \in \R} \left\{
                \alpha + \frac1\theta \left(e ^{\theta (x - \alpha)} - 1\right) :|\hX - x| \le \rho
            \right\}
        \right] \\
        &= \inf _{\alpha \in \R} \E _{\hX \sim \hP} \left[ 
            \alpha + \frac1\theta \left(e ^{\theta (\hX + \rho - \alpha)} - 1\right) - \alpha
        \right] \\
        &= \frac1\theta\log (\E _{\hX \sim \hP} [e ^{\theta \hX}]) + \rho = \EnthP _\theta (\hX) + \rho.
    \end{align*}

    For $p < \infty$, we verify $\cLJ (\rho) = +\infty$ directly. We define $\P _\epsilon = (1 - \epsilon) \hP + \epsilon \hP _{M _\epsilon}$, where $\hP _M = (x \mapsto x + M) _\# \hP$ is right-translation of $\hP$ by $M$, and $M _\epsilon = \rho \epsilon ^{-1/p}$. Then $\cW _p (\hP, \P) \le \rho$. However, 
    \begin{align*}
        \E _{\P _\epsilon} [e ^{\theta X}] = (1 - \epsilon) \E [e ^{\theta \hX}] + \epsilon \E [e ^{\theta (\hX + M _\epsilon)}] = (1 - \epsilon + \epsilon e ^{\theta M _\epsilon}) \E [e ^{\theta \hX}], 
    \end{align*}
    and accordingly 
    \begin{align*}
        \Ent ^{\P _\epsilon} _\theta (X) = \frac1\theta\log (\E _{X \sim \P _\epsilon} [e ^{\theta X}]) = \frac1\theta\log (\E _{\hX \sim \hP} [e ^{\theta \hX}]) + \frac1\theta \log (1 + \epsilon (e ^{\theta \epsilon ^{-1/p}} - 1)),
    \end{align*}
    which tends to infinity as $\epsilon \to 0$.
    \QED
\endproof

\section{Proofs for Section \ref{sec:globalized}}
\label{app:globalized}

\proof{Proof of Proposition \ref{prop:globalized}.}

For a fixed $\theta$, first we apply Theorem \ref{thm} to $-\cLG (\cdot, \theta)$ by taking the Fenchel conjugate
\begin{align*}
    (-\cLG (\cdot, \theta)) ^* (-\lambda) &= \sup _{\P, \tP \in \cP (\cX)} \left\{ 
        \E _{X \sim \P} [f (X)] - \lambda \cK (\tP, \P): \tcK (\hP, \tP) \le \theta
    \right\} \\
    &= \sup _{\tP \in \cP (\cX)} \left\{ 
        \E _{\tX \sim \tP} \left[
            \sup _{x \in \cX} \left\{ f (x) - \lambda c (\tX, x)
            \right\}
        \right]: \tcK (\hP, \tP) \le \theta
    \right\}.
\end{align*}
Denote $\tf (\tx) = \sup _{x \in \cX} f (x) - \lambda c (\tx, x)$. Then we apply Theorem \ref{thm} to $\tf, \tc$ and $-(-\cLG (\cdot, \theta)) ^* (-\lambda)$ by taking Fenchel conjugate $\theta \to -\mu$:
\begin{align*}
    \cLG \doublestar (-\lambda, -\mu) &= \sup _{\tP \in \cP (\cX)} \left\{ 
        \E _{\tX \sim \P} \left[
            \sup _{x \in \cX} \left\{
                f (x) - \lambda c (\tX, x)
            \right\}
        \right] - \mu \tcK (\hP, \tP)
    \right\} \\
    &= \E _{\hX \sim \hP} \left[
            \sup _{x, \tx \in \cX} \left\{
                f (x) - \lambda c (\tx, x) - \mu \tc (\hX, \tx)
            \right\}
        \right].
\end{align*}
Since $(-\cLG (\cdot, \theta)) ^* (-\lambda)$ is concave in $\theta$, we recover it by 
\begin{align*}
    (-\cLG (\cdot, \theta)) ^* (-\lambda) &= \min _{\mu \ge 0} \left\{
        \mu \theta + (-\cLG) \doublestar (-\lambda, -\mu)
    \right\} \\
    &= \min _{\mu \ge 0} \left\{
        \mu \theta + \E _{\hX \sim \hP} \left[
            \sup _{x, \tx \in \cX} \left\{
                f (x) - \lambda c (\tx, x) - \mu \tc (\hX, \tx)
            \right\}
        \right]
    \right\}.
\end{align*}
Since $\cLG (\rho, \theta)$ is concave in $\rho$, we recover it by 
\begin{align*}
    \cLG (\rho, \theta) = \min _{\lambda, \mu \ge 0} \left\{
        \lambda \rho + \mu \theta + \E _{\hX \sim \hP} \left[
            \sup _{x, \tx \in \cX} \left\{
                f (x) - \lambda c (\tx, x) - \mu \tc (\hX, \tx)
            \right\}
        \right]
    \right\}.
\end{align*}
In particular, if $c (x _1, x _2) = \tc (x _1, x _2) = d (x _1, x _2)$ are the same metric, then 
\begin{align*}
    \cLG \doublestar (-\lambda, -\mu) &= \E _{\hX \sim \hP} \left[
        \sup _{x, \tx \in \cX} \left\{
            f (x) - \lambda d (\tx, x) - \mu d (\hX, \tx)
        \right\}
    \right] \\
    &= \E _{\hX \sim \hP} \left[
        \sup _{x \in \cX} \left\{
            f (x) - (\lambda \wedge \mu) d (\hX, x)
        \right\}
    \right]
\end{align*}
by taking $\tx = x$ when $\lambda \ge \mu$ and $\tx = \hX$ when $\lambda \le \mu$. Correspondingly, 
\begin{align*}
    (-\cLG) ^* (-\lambda, \theta) 
    &= \min _{\mu \ge 0} \left\{
        \mu \theta + \E _{\hX \sim \hP} \left[
            \sup _{x \in \cX} \left\{
                f (x) - (\lambda \wedge \mu) d (\hX, x)
            \right\}
        \right]
    \right\} \\
    &= \min _{\mu \in [0, \lambda]} \left\{
        \mu \theta + \E _{\hX \sim \hP} \left[
            \sup _{x \in \cX} \left\{
                f (x) - \mu d (\hX, x)
            \right\}
        \right]
    \right\},
\end{align*}
and 
\begin{align*}
    \cLG (\rho, \theta) &= \min _{0 \le \mu \le \lambda} \left\{
        \lambda \rho + \mu \theta + \E _{\hX \sim \hP} \left[
            \sup _{x \in \cX} \left\{
                f (x) - \mu d (\hX, x)
            \right\}
        \right]
    \right\} \\
    &= \min _{\mu \ge 0} \left\{
        \mu (\rho + \theta) + \E _{\hX \sim \hP} \left[
            \sup _{x \in \cX} \left\{
                f (x) - \mu d (\hX, x)
            \right\}
        \right]
    \right\}.
    \tag*{\QED}
\end{align*}
\endproof

\end{document}